\documentclass[a4paper]{article}
\usepackage{amsmath}
\usepackage{graphics}
\usepackage{latexsym}
\usepackage{amsfonts}
\usepackage{mathrsfs}
\usepackage{amssymb}
\usepackage{amssymb,amsmath,amsthm, amsfonts}
\numberwithin{equation}{section}

\newtheorem{propo}{Proposition}[section]

\newtheorem{lemma}{Lemma}[section]

\newtheorem{assumption}{Assumption}
\def \u {\mathcal{U}}
\def \rwk {\rightarrow_k}
\def\qed{\ \vrule width.2cm height.2cm depth0cm\smallskip}
\def \et{\eta}

\def \x{X^{t,x}}
\def \ed {\end{document}}
\def \tx {(t,x)\in \esp}

\def \fl {\forall}
\def \lb{\label}
\def\wr {w.r.t.}

\def \et {\eta}

\def \ie {\int_E}

\newcommand{\eps}{\varepsilon}

\newcommand{\brm}{\begin{rem}}
\newcommand{\ermq}{\end{rem}}

\newcommand{\ba}{\begin{array}}
\newcommand{\ea}{\end{array}}
\newcommand{\be}{\begin{equation}}
\newcommand{\ee}{\end{equation}}
\newcommand{\bea}{\begin{eqnarray}}
\newcommand{\eea}{\end{eqnarray}}
\newcommand{\beaa}{\begin{eqnarray*}}
\newcommand{\eeaa}{\end{eqnarray*}}

\def \R{I\!\!R}
\def \E{\mathbb{E}}

\def\b{\beta}
\def\g{\gamma}
\def\d{\delta}

\def\z{\zeta}

\def\l{\lambda}

\def\si{\sigma}
\def\t{\tau}

\def\D{\Delta}

\def\cA{{\cal A}}
\def\cB{{\cal B}}
\def\cC{{\cal C}}

\def\cE{{\cal E}}
\def\cF{{\cal F}}

\def\cH{{\cal H}}

\def\cL{{\cal L}}

\def\cP{{\cal P}}

\def\cS{{\cal S}}

\def\cU{{\cal U}}

\def\no{\noindent}

\def\ms{\medskip}

\def\q{\quad}
\def\qq{\qquad}

\def \bp {\bold P}
\def \ife {(1\wedge |e|^2)_{e\in E}}
\def\qed{ \hfill \vrule width.25cm height.25cm depth0cm\smallskip}

\newcommand{\basa}{\begin{assumption}}
\newcommand{\easa}{\end{assumption}}

\newcommand{\bas}{\begin{assum}}
\newcommand{\eas}{\end{assum}}

\def \xk {^k\!\!X^{t,x}}
\def \yk {^k\!Y^{t,x}}
\def \zk {^k\!\!Z^{t,x}}
\def \uk {^k\!U^{t,x}}
\def \zki {^k\!\!Z^{i;t,x}}

\def \xkk {^k\!\!X^{t,x_k}}

\def\dis{\displaystyle}

\def \P{\mathbb{P}}
\newtheorem{thm}{Theorem}[section]

\newtheorem{prop}[thm]{Proposition}
\newtheorem{rem}[thm]{Remark}

\newtheorem{defn}[thm]{Definition}
\newtheorem{assum}[thm]{Assumption}

\def\mb{\mbox}
\def \im {{i=1,m}}
\def \rw {\rightarrow}
\def \itk {1_{\{|e|\ge \frac{1}{k}\}}}
\def \R{\mathbb{R}}

\def \esp {[0,T]\times \R^k}

\def \r {\mathbb{R}}
\def \rk {\r^k}
\def \mo {\{1,\dots,m\}}

\def \pg {\Pi_g}

\def \pgc {\Pi_g^c}
\def \yx {Y^{t,x}}
\def \zx {Z^{t,x}}
\def \ux {U^{t,x}}

\title{Viscosity solutions of second order integral-partial differential equations without monotonicity condition: A new result} 
\author{
Said Hamad\`ene\thanks{Universit\'e du
Maine, LMM, Avenue Olivier Messiaen, 72085 Le Mans, Cedex 9, France, e-mail: hamadene@univ-lemans.fr}  }

\begin{document}
\date{\today}
\maketitle

\begin{abstract}
We show existence and uniqueness of a continuous with polynomial growth viscosity solution of a system of second order integral-partial differential equations (IPDEs for short) without assuming the usual monotonicity condition of the generator with respect to the jump component as in Barles et al.'s article \cite{BarlesBuckPardoux}. The L\'evy measure is arbitrary and not necessarily finite. In our study the main tool we used is the notion of backward stochastic differential equations with jumps. 
 \end{abstract}
\no{\bf AMS Classification subjects}: 35D40 ; 35K10 ;
60H30.
\medskip

\no {$\bf Keywords$}: Integral-partial differential equation ; Backward stochastic differential equation with jumps ; Viscosity solution ; Non-local operator.
\section{Introduction}
The main objective of this paper is to deal with the following system of integral-partial differential equations: $\forall i\in \mo$, 
\be \left\{ \label{ipde-intro} \begin{array}{l} 
-\partial_t u^i(t,x)- b(t,x)^\top D_xu^i(t,x)-\frac{1}{2}\mathrm{Tr}\big(\sigma\sigma^\top (t,x) D^2_{xx}u^i(t,x)\big)-Ku^i(t,x)\\\q\qq -h^{(i)}(t,x, (u^j(t,x))_{j=1,m}, (\sigma^\top D_xu^i)(t,x), B_iu^i(t,x))=0,\,\, (t,x)\in \esp ;\\\\
u^i(T,x) = g^i(x)
\end{array}  \right.
\ee where the operators $B_i$ and $K$ are defined as follows:\be\label{defbk-intro}\begin{array}{l} B_iu^i(t,x)=\int_E\g^i(t,x,e)\left( u^i(t,x
+\beta(t,x,e))-u^i(t,x)\right)\l(de)\mbox{ and }\\\\
Ku^i(t,x)=\int_E\left( u^i(t,x
+\beta(t,x,e))-u^i(t,x)-\beta(t,x,e)^\top 
D_xu^i(t,x)\right)\l(de)\end{array}\ee  
where $\l$ is a L\'evy measure on $E:=\R^\ell -\{0\}$ which integrates the function $(1\wedge |e|^2)_{e\in E}$. 

The second order system of equations (\ref{ipde-intro}) is of non-local type since the operators $B_iu^i$ and $Ku^i$ at $(t,x)$ involve the values of $u_i$ in the whole space $\r^k$ and not only locally, i.e. in a neighbourhood of $(t,x)$. 

This system  of IPDEs, introduced by Barles et al. in \cite{BarlesBuckPardoux}, is deeply related to the following multidimensional backward stochastic differential equation (BSDE for short) with jumps whose solution, for fixed $\tx$, is a triple of adapted stochastic processes $(\yx_s,\zx_s,\ux_s)_{s\le T}$ with values in $\r^m\times \r^{m\times d}\times L^2(\l)$ which mainly satisfy:  
$\forall i\in \mo$, 
\be \label{bsde-jumps-intro(o)}
   \left\{ \begin{array}{l} 
-dY_s^{i;t,x} = h^{(i)}(s,X^{t,x}_s,(Y^{j;t,x}_s)_{j=1,m},Z_s^{i;t,x}, \ie \g_i(s,\x_s,e)U_s^{i;t,x}(e)\l(de))ds \\\qq\qq \qq -Z_s^{i;t,x} dB_s -\int_{E} 
    U_{s}^{i;t,x}(e)\tilde{\mu}(ds,de),\,\,\forall \,s\leq T\,;\\\\Y_T^{i;t,x}= g^{i}(X^{t,x}_T),\end{array}\right.\ee
where: 

(i) $B:=(B_s)_{s\leq T}$ is a $d$-dimensional Brownian motion, $\mu$ an independent Poisson random measure 
with compensator $ds\l(de)$ and $\tilde \mu (ds, de):=\mu (ds, de)-ds\l(de)$ its compensated random measure ; 

(ii) 
for any $(t,x)\in \esp$, $(X^{t,x}_s)_{s\leq T}$ is the solution of the following standard stochastic differential equation of diffusion-jump type, i.e., 
\be \label{eq:diffusion(o)}\begin{array}{l}X_{s}^{t,x} = x + \int_{t}^{s} b(r,X^{t,x}_r)dr +\int_{t}^{s}\sigma(r, X^{t,x}_r)dB_r
 +\int_{t}^{s}\int_{E}\beta(r, X^{t,x}_{r-}, e) \tilde{\mu}(dr,de), \mbox{ for } s\in [t,T] \mbox{ and }X_s^{t,x}=x \mbox{ if }s\leq t. \end{array}\ee
Actually it has been shown in \cite{BarlesBuckPardoux} that, under standard assumptions on the functions $b$, $\sigma$, $\beta$, $g^i$, $h^{(i)}$ and $\g_i$ and due to the Markovian framework of randomness which stems from the Markov process $X^{t,x}$, there exist deterministic continuous functions $(u^i(t,x))_{i=1,m}$ such that for any $s\in [t,T]$, \be \label{fk}Y^{i;t,x}_s=u^i(s,\x_s), \forall i=1,...,m.\ee Moreover if for any $i=1,...,m$, 

(a) $\g_i\ge 0$ ;

(b) the mapping $q\in \R\mapsto h^{(i)}(t,x,y,z,q)$ is non-decreasing, when the other components $(t,x,y,z)$ are fixed ;  

\noindent then the functions $(u^i)_{i=1,m}$ is the unique continuous viscosity solution of system (\ref{ipde-intro}) in the class of functions with polynomial growth (at least). Conditions (a)-(b), which will be referred as the monotonicity conditions, are needed in \cite{BarlesBuckPardoux} in order to have the comparison property and to treat the operator $B_iu^i$ which is not well-defined for an arbitrary $u$. However we should point out those conditions are not required in order to show the existence and uniqueness of the solution $(\yx,\zx,\ux)$ of BSDE (\ref{bsde-jumps-intro(o)}). 

Therefore the main issue is to deal with the viscosity solutions of system (\ref{ipde-intro}) without assuming the above conditions (a)-(b) neither on $\g_i$ nor on $h^{(i)}$, $i=1,...,m$. A step forward in the resolution of this problem is made by Hamad\`ene-Morlais in \cite{Hammorlais-mai2015} where it is shown that, when the L\'evy measure $\l$ is finite i.e. $\l(E)<\infty$, then system (\ref{ipde-intro}) has a unique solution which is given by the functions $(u^i)_{i=1,m}$ defined in (\ref{fk}).  

The main objective of this paper is once more to deal with the problem of existence and uniqueness of a viscosity solution of system of IPDEs (\ref{ipde-intro}) without assuming the monotonicity conditions neither on $\g_i$ nor on $h^{(i)}$, $i=1,...,m$ and for an arbitrary L\'evy measure $\l$ without assuming its finitness as in \cite{Hammorlais-mai2015}. There are two crucial points. The first one is the characterization (\ref{expressiondessauts}) below of the process $\ux=(U^{i;t,x})_{i=1,m}$ of the solution of the BSDE (\ref{bsde-jumps-intro(o)}) by means of the functions $(u^i)_{i=1,m}$ defined in (\ref{fk}) and the jump-diffusion process $X^{t,x}$. Actually, using the truncation method at the origin of the L\'evy measure $\l$ we show that  for any $i=1,...,m$, 
\be \lb{expressiondessauts}
U^{i;t,x}_s(e)=u^i(s,X^{t,x}_{s-}+\beta(s,X^{t,x}_{s-},e))-u^i(s,X^{t,x}_{s-}),\,\,ds\otimes d\P\otimes d\l \mbox{ on }[t,T]\times \Omega \times E.
\ee
The second one is the local boundedness  of the increment rate w.r.t $x$ of the functions $u^i$ which is obtained under reasonnable conditions on the functions $h^{(i)}$ and $\gamma_i$. Those facts allow us to avoid to replace $B_iu^i$ with $B_i\phi$ where $\phi$ is the test function, as in \cite{BarlesBuckPardoux}. We then introduce a new definiton of the viscosity solution of system (\ref{ipde-intro}) and relying on Barles et al.'s result \cite{BarlesBuckPardoux} and, on the other hand, on BSDEs with jumps ones we show that the functions defined in (\ref{fk}) is the unique viscosity solution of system (\ref{ipde-intro}). 
Our definition of a viscosity solution of (\ref{ipde-intro}) is not the same as the one in \cite{BarlesBuckPardoux} and looks like to the one given in \cite{Hammorlais-mai2015}. This is the novelty of this paper and according to our best knowledge this result is not obtained yet in a so general framework.  
\ms

Note that there are also other papers on this topic of IPDEs amongst one can quote (\cite{Alvareztourin, BarlesChassImb08, BarlesImbert08, bucdahnhuli}, etc. and the references therein).  Finally let us point out that IPDEs which do not satisfy the monotonicity conditions are encountered in mathematical finance when dealing with the problem of liquidation of portfolios (see e.g. \cite{horst}).  
\ms

This paper is organized as follows. Section 2 is devoted to fix the framework on which we are working and, for completeness, to recall the state of the art on the main subject. Section 3 is mainly devoted to the proof of the relation (\ref{expressiondessauts}). We first prove that the increment rates of the functions $u^i$, $i=1,...,m$, are locally bounded. Later on, by the method of truncation of the L\'evy measure $\l$ at the origin in such a way to get into the setting of a finite L\'evy measure which is already considered in \cite{Hammorlais-mai2015}, we prove by approximations the relation (\ref{expressiondessauts}). In Section 4 we precise the notion of viscosity solution we are working with and we give the proof of the main result. We emphasize that this definition is not the same as the one in \cite{BarlesBuckPardoux}. Finally new types of systems of IPDEs are introduced and discussed in Section 5. \qed

\section{Framework and state of the art}
Let $(\Omega,\mathcal{F},(\mathcal{F}_t)_{t\leq T}, \mathbb{P})$ be a stochastic basis such that $\cF_0$ contains all $\P$-null sets of $\cF$, and $\cF_{t}=\cF_{t+}:=\bigcap_{\eps> 0}\cF_{t+\eps}$, $t\geq 0$, and we suppose that the filtration is generated by the two mutually independant processes:

\no (i) $B:=(B_t)_{t\ge 0}$ a $d$-dimensional Brownian motion and

\no (ii) a Poisson random measure $\mu$ on $\r^+\times E$, where $E:=\r^{\ell}-\{0\}$ is equipped with its Borel field $\cE$ ($\ell \ge 1$).
The compensator $\nu(dt,de)=dt\l(de)$ is such that $\{\tilde \mu([0,t]\times A)=(\mu-\nu)([0,t]\times A)\}_{t\geq 0}$ 
is a martingale for all $A\in \cE$ 
satisfying $\l(A)<\infty$. We also assume that $\l$ is a $\sigma$-finite measure on $(E,\cE)$, integrates the function $(1\wedge |e|^2)_{e\in E}$ and $\l(E) = \infty$. Note that the case when $\l(E)<\infty$ is already considered in \cite{Hammorlais-mai2015}. 
\ms

Next we denote by: 

\no (iii) $\cP$ (resp. $\bp$) the field on $[0,T]\times \Omega$ 
of $(\cF_t)_{t\leq T}$-progressively measurable (resp. predictable) sets ;

\no (iv) For $\kappa \ge 1$, $L^2_\kappa(\l)$ the space of Borel measurable functions $\varphi:=(\varphi(e))_{e\in E}$ from $E$ into $\r^\kappa$ such that $\|\varphi\|_{L^2_\kappa (\l)}^2:=\int_E|\varphi(e)|_\kappa^2\l(de)<\infty$ ; $L^2_1(\l)$ will be simply denoted by $L^2(\l)$ ; 

\no (v) $\mathcal{S}^{2}(\mathbb{R}^{\kappa})$ the space of RCLL (for right continuous with left limits) $\cP$-measurable and $\r^\kappa$-valued processes such that
$\mathbb{E}[\sup_{s\le T} |Y_s|^2] <\infty $ ; $\cA^2_c$ is its subspace of continuous non-decreasing processes $(K_t)_{t\leq T}$ such that $K_0=0$ ;  

\no (vi) $\mathcal{H}^2(\mathbb{R}^{\kappa\times d})$ the space of processes $Z:=(Z_s)_{s\le T}$ which are $\cP$-measurable, $\r^{\kappa \times d}$-valued and satisfying \\$\mathbb{E}[ \int_{0}^{T} |Z_s|^{2} ds ]<\infty $ ; 

\no (vii) 
$\mathcal{H}^2(L_\kappa^2(\l))$ the space of processes $U:=(U_s)_{s\le T}$ which are $\bp$-measurable, $L_\kappa^2(\l)$-valued and satisfying \\
$\mathbb{E}[ \int_{0}^{T} \|U_s(\omega)\|_{L_\kappa^2(\l)}^{2} ds ]<\infty$ ; 

\no (viii) $\pg$ the set of deterministic functions $\varpi$: $(t,x)\in \esp \mapsto \varpi(t,x)\in \r$ of polynomial growth, i.e., for which there exists two non-negative constants $C$ and $p$ such that for any $(t,x)\in \esp$, 
$$
|\varpi(t,x)|\leq C(1+|x|^p).
$$The subspace of $\pg$ of continuous functions will be denoted by $\pgc$ ;

\no (ix)  $\mathcal{U} $ the subclass of $\Pi_g^c$ which consists of functions $\Phi: (t,x)\in \esp \longmapsto \r$ such that for some non-negative constants $C$ and $p$ we have
$$|\Phi(t, x) - \Phi(t,x')| \le C(1+|x|^p+|x'|^p)|x-x'|,\,\, \mb{for any }t,x,x'. 
$$

\no (x) For any process $\theta:=(\theta_s)_{s\le T}$ and $t\in (0,T]$, $\theta_{t-}=\lim_{s\nearrow t}\theta_s$ and 
$\D_t \theta=\theta_t-\theta_{t-}$ ;
\ms

Now let $b$ and $\sigma$ be the following functions:
$$\ba{l}
b:(t,x)\in \esp \mapsto b(t,x) \in \r^k\\
\si :(t,x)\in\esp \mapsto \si (t,x)\in \r^{k\times d}.\ea
$$
We assume that they are jointly continuous in $(t,x)$ and Lipschitz continuous $\wr$ $x$ uniformly in $t$, i.e., there exists a constant $C$ such that  
\begin{equation}\label{bslip}
\forall \; (t,x, x') \in  [0,T] \times \mathbb{R}^{k+k},\;\; |b(t,x) -b(t,x'| +|\sigma(t,x) -\sigma(t, x')| \le C |x -x'|.
\end{equation}
Let us notice that by (\ref{bslip}) and continuity, the functions $b$ and $\si$ are of linear growth, i.e., there exists a constant $C$ such that  
\begin{equation}\label{bslip2}
\forall \; (t,x) \in [0,T] \times \mathbb{R}^d \;\;|b(t,x)| +|\sigma(t,x)|\le C(1+|x|).
\end{equation}
Let $\beta:(t,x,e)\in \esp \times E \mapsto \beta (t,x,e)\in \rk$ be a measurable function such that for some real constant $C$, and for all $e\in E$, 
\begin{equation}\label{cdbeta}
\begin{array}{l}
\mbox{(i)} \q |\beta(t,x,e)|\le C (1\wedge |e|);\\
\mbox{(ii)}\q |\beta(t,x,e)-\beta (t,x',e)|\le C |x-x'|(1\wedge |e|);\\
\mbox{(iii)}\q \mbox{the mapping }(t,x)\in \esp  \rw \beta (t,x,e)\in \rk \mbox{ is continuous for any } e\in E.
\end{array}
\end{equation}

\no Once for all, throughout this paper, we assume that conditions (\ref{bslip}), (\ref{bslip2}) and (\ref{cdbeta}), on $b$, $\sigma$ and $\beta$ respectively, are fulfilled. \ms

Next let $(t,x)\in \esp$ and $(X^{t,x}_s)_{s\leq T}$ be the stochastic process solution of the following standard stochastic differential equation of diffusion-jump type:
\be \label{eq:diffusion}\begin{array}{l}X_{s}^{t,x} = x + \int_{t}^{s} b(r,X^{t,x}_r)dr +\int_{t}^{s}\sigma(r, X^{t,x}_r)dB_r
 +\int_{t}^{s}\int_{E}\beta(r, X^{t,x}_{r-}, e) \tilde{\mu}(dr,de), \mbox{ for } s\in [t,T] \mbox{ and }X_s^{t,x}=x \mbox{ if }s\leq t. \end{array}\ee
Under assumptions (\ref{bslip}), (\ref{bslip2}) and (\ref{cdbeta}) the solution of equation (\ref{eq:diffusion}) exists and is unique (see \cite{kunita} for more details). Moreover it satisfies the following estimates: $\forall p\geq 2$, $x,x'\in \r^k$ and $s\ge t$,
\be \label{estimx}
\E [\sup_{r\in [t,s]}|X^{t,x}_r-x|^p]\leq M_p(s-t)(1+|x|^p)] \mbox{ and }  
\E [\sup_{r\in [t,s]}|X^{t,x}_r-X^{t,x'}_r-(x-x')|^p]\leq M_p(s-t)|x-x'|^p
\ee
for some constant $M_p$. $\Box$
\ms

We are now going to introduce the objects which are  specifically connected to the BSDEs with jumps we will deal with. Let $(g^{i})_{i=1,m}$ and $(h^{(i)})_{i=1,m}$ be functions defined as follows: For $i=1,...,m$, 
$$\ba{ccc}\ba{cl}
g^i: \r^k &\longrightarrow \r^m\\
\qq x&\longmapsto g^i(x) \ea &\mb{and}\qq \ba{cl}
h^{(i)}: [0,T]\times \r^{k+m+d+1}& \longrightarrow \r \\
(t,x,y,z,q)&\longmapsto h^{(i)}(t,x,y,z,q).
\ea\ea
$$
Moreover we assume they satisfy:
\ms

\noindent {\bf (H1)}: For any $i\in \mo$, the function $g^i$ belongs to $\u$. 
\ms

\noindent {\bf (H2)}: For any $i\in \mo$, 

\no (i) the function $h^{(i)}$ is Lipschitz in $(y,z,q)$ uniformly in $(t,x)$, i.e., there exists a real constant $C$ such that for any $(t,x)\in \esp$, $(y,z,q)$ and $(y',z',q')$ elements of $\r^{m+d+1}$, 
\be \label{cdlipschitzfyzu}   |h^{(i)}(t,x, y, z,q) -h^{(i)}(t,x, y', z',q')|
\leq C(|y-y'|+|z-z'|+|q-q'|); \ee 

\no (ii) the function $(t,x)\mapsto h^{(i)}(t,x,y,z,q)$, for fixed $(y,z,q)\in \r^{m+d+1}$, belong uniformly to $\u$, i.e., it is continuous and there exist constants $C$ and $p$ (which do not depend on $(y,z,q)$) such that, \be \lb{polygrtf}|h^{(i)}(t, x,y,z,q) - h^{(i)}(t, x',y,z,q)| \le C(1+|x|^p+|x'|^p)|x-x'|,\,\, \mb{for any }t,x,x'.
\ee

\no Next let $\g_i$, $i=1,\dots, m$ be Borel measurable functions defined from $\esp \times E$ into $\r$ and satisfying: 
\be \label{condgamma}
\begin{array}{l}
(i)\q |\gamma_i(t,x,e)|\leq C(1\wedge |e|)\\
(ii)\q |\gamma_i(t,x,e)-\gamma_i(t,x',e)|\leq 
C(1\wedge |e|)|x-x'|(1+|x|^p+|x'|^p)
\\
(iii) \mbox{ the mapping } t\in [0,T]\longmapsto \g_i(t,x,e) \; \mbox{is continuous for any } (x,e). \end{array}
\ee          

\no Finally let us introduce the following functions 
$(f^{(i)})_{i=1,m}$, defined by:
\be\lb{deffi}\ba{l}
\forall (t,x,y,z,\z)\in [0,T]\times \r^{k+m+d}\times L^2(\l), \,\,
f^{(i)}(t,x,y,z,\z):=h^{(i)}(t,x,y,z,\ie \g_i(t,x,e)\z(e)\l(de)).\ea
\ee
The functions $f^{(i)}$, $i=1,\dots, m$, enjoy the two following properties: 
\ms 
\be \lb{condfi}\ba{l}
\mb{(a) } f^{(i)}\mb{ is Lipschitz in }(y,z,\z), \mb{ uniformly in }(t,x),\mb{ i.e., there exists a constant }C\mb{ such that } \\\\
 \qq\qq |f^{(i)}(t,x, y, z,\z) -f^{(i)}(t,x, y', z',z')|
\leq C(|y-y'|+|z-z'|+\|\z-\z'\|_{L^2(\l)})\\\\
\mb{ since }h^{(i)}\mb{ is uniformly Lipschitz in }(y,z,q)\mb{ and }\g_i\mb{ verifies } (\ref{condgamma})-(i).\\\\

\mb{ (b) the function }(t,x)\in \esp \mapsto f^{(i)}(t,x,0,0,0)\mb{ belongs to }\pg^c\mb{ and then }\E[\int_0^T|f^{(i)}(r,\x_r,0,0,0)|^2dr]<\infty.  
\ea
\ee

Let now $(t,x)\in \esp$ and let us consider the following $m$-dimensional BSDE with jumps: 
\be \left\{ \label{mainbsde}
    \begin{array}{l}  
\vec{Y}^{t,x}:=(Y^{i;t,x})_{i=1,m}\in \cS^2(\r^m),\,
{Z}^{t,x}:=({Z}^{i;t,x})_{i=1,m}\in \cH^2(\r^{m\times d}), U^{t,x}:=(U^{i;t,x})_{i=1,m}\in \cH^2(L_m^2(\l));\\ \forall i\in \mo, \,
Y_T^{i}= g^{i}(X^{t,x}_T) \mbox{ and }\\
\qq dY_s^{i;t,x} = -f^{(i)}(s,X^{t,x}_s, \vec{Y}^{t,x}_s,Z_s^{i;t,x}, U_s^{i;t,x})ds -Z_s^{i;t,x} dB_s -\int_{E} 
    U_{s}^{i;t,x}(e)\tilde{\mu}(ds,de),\,\,\forall s\leq T,\end{array}\right.  \ee
 where for any $i\in \mo$, $Z_s^{i;t,x}$ is the $i$-th row of $Z_s^{t,x}$ and $U_s^{i;t,x}$ is the $i$-th component of $U_s^{t,x}$.
 \ms
 
The following result is related to existence and uniqueness of a solution for the BSDE with jumps (\ref{mainbsde}). Its proof is given in 
Li-Tang \cite{TangLi94} (one can also see Barles et al. \cite{BarlesBuckPardoux}).
\ms

\begin{propo} \label{existencegene}(Tang-Li, \cite{TangLi94}): Assume that Assumptions (H1)-(H2) hold. Then for any $(t,x)\in \esp$,
the BSDE (\ref{mainbsde}) has a unique solution 
$(\vec{Y}^{t,x},{Z}^{t,x},U^{t,x})$. \end{propo}

\begin{rem}\lb{remexistuniq} The solution of this BSDE exists and is unique since:

(i) $\E[|g(\x_T)|^2]<\infty$, due to polynomial growth of $g$ and estimate (\ref{estimx}) on $\x$ ;

(ii) for any $i=1,\dots,m$,  $f^{(i)}$ verifies the properties (\ref{condfi})-(a),(b) related to uniform Lipschitz w.r.t $(y,z,\z)$ and $ds\otimes d\P$-square integrability of the process $(f^{(i)}(s,\x_s,0,0,0))_{s\le T}$. \qed 
\end{rem}

Next, the following result proved in Barles et al. (\cite{BarlesBuckPardoux}, Proposition 2.5 and 
Theorems 3.4, 3.5), establishes the relationship between the solution of (\ref{mainbsde}) 
and the one of system (\ref{ipde-intro}).

\begin{propo}\label{solvisco}(\cite{BarlesBuckPardoux}): Assume that (H1) and (H2) are fulfilled. Then there exist deterministic continuous functions $(u^i(t,x))_{i=1,m}$ which belongs to $\pg$ such that for any $(t,x)\in \esp$, the solution of the BSDE (\ref{mainbsde}) verifies: 
\be\label{repre}
\forall i\in \mo,\,\,\forall s\in [t,T],\,\,Y^{i;t,x}_s=u^i(s,X^{t,x}_s).
\ee
Moreover if for any $i\in \mo$, 

(i) $\g^i\geq 0$ ;

(ii) for any fixed $(t,x,\vec{y},z)\in [0,T]\times \r^{k+m+d}$, the mapping
$
q\in \r\longmapsto h^{(i)}(t,x,\vec{y},z,q)\in \r
$
is non-decreasing. \\ Then the functions $(u^i)_{i=1,m}$ is a continuous viscosity solution (in Barles et al.'s sense, see Definition \ref{bbpdef} in Appendix) of the following system of IPDEs: $\forall i\in \mo$, 
\be \left\{ \label{secondorder-pde} \begin{array}{l} 
-\partial_t u^i(t,x)- b(t,x)^\top D_xu^i(t,x)-\frac{1}{2}\mathrm{Tr}\big(\sigma\sigma^\top (t,x) D^2_{xx}u^i(t,x)\big)-Ku^i(t,x)\\\q\qq -h^{(i)}(t,x, (u^j(t,x))_{j=1,m}, (\sigma^\top D_xu^i)(t,x), B_iu^i(t,x))=0,\,\, (t,x)\in \esp ;\\
u^i(T,x) = g^i(x),
\end{array}  \right.
\ee where \be\label{defbk}\begin{array}{l} B_iu^i(t,x)=\int_E\g^i(t,x,e)\{u^i(t,x
+\beta(t,x,e))-u^i(t,x)\}\l(de)\mbox{ and }\\\\
Ku^i(t,x)=\int_E\{u^i(t,x
+\beta(t,x,e))-u^i(t,x)-\beta(t,x,e)^\top 
D_xu^i(t,x)\}\l(de).\end{array}\ee 
Finally, the solution $(u^i(t,x))_{i=1,m}$ of 
(\ref{secondorder-pde}) is unique in the class $\Pi^c_g$. \qed
\end{propo}
\begin{rem}
(i) The solution $u=(u^i)_{i=1,m}$ is also unique in the class of functions which satisfy the following weaker growth condition:
$$
\lim_{|x|\rw \infty}|u(t,x)|e^{-\tilde A [ln(|x|)]^2}=0$$
uniformly for $t\in [0,T]$, for some $\tilde A > 0$ (see \cite{BarlesBuckPardoux} or \cite{BarlesImbert08} for more details).\\
(ii) The functions $h^{(i)}$ verify the condition (A2.v) in (\cite{BarlesBuckPardoux}, pp. 73), 
under which uniqueness of the solution of (\ref{secondorder-pde}) is obtained, by the assumption (H2)-(ii). 
\end{rem}
\section{Estimates and properties}
Our next step is to provide estimates for the functions $(u^i)_{i=1,m}$ defined in (\ref{repre}). Recall that $(\vec{Y}^{t,x},\zx,\ux):=((Y^{i;t,x})_{i=1,m},(Z^{i;t,x})_{i=1,m},(U^{i;t,x})_{i=1,m})$ is the unique solution of the BSDE with jumps (\ref{mainbsde}).
\begin{lemma}\label{uplem}
Under (H1)-(H2), for any $p\geq 2$ there exist two non-negative constants $C$ and $\rho$ such that
\be\label{up}
\E\Big [\Big\{\int_0^Tds(\int_E|U^{t,x}_s(e)|^2\l (de))\Big \}^{\frac{p}{2}}\Big ]=
\E\Big [\Big\{\int_0^Tds\|U_s^{t,x}\|^2_{L^2_m(\l)}\Big \}^{\frac{p}{2}}\Big ]\leq C(1+|x|^\rho).\ee
\end{lemma}
\proof First let us point out that since $X^{t,x}_s=x$ for $s\in [0,t]$ then, uniqueness of the solution of BSDE (\ref{mainbsde}) implies that 
\be\label{uzavant}
Z^{t,x}_s=0 \mbox{ and }U^{t,x}_s=0,\,\,ds\otimes d\bp-a.e. \,\mbox{ on }[0,t]\times \Omega. \ee
Next let $p\geq 2$ be fixed. Using the representation (\ref{repre}), for any $i\in \mo$ and $s\in [t,T]$ we have 
\be\label{bsdeinter}
Y^{i;t,x}_s=g^i(X^{t,x}_T)+
 \int_s^Tf^{(i)}(r,X^{t,x}_r, (u^j(r,X^{t,x}_r))_{j=1,m},Z_r^{i;t,x}, U_r^{i;t,x})dr -\int_s^TZ_r^{i;t,x} dB_r-\int_s^T\int_{E} 
    U_{r}^{i;t,x}(e)\tilde{\mu}(dr,de).
\ee
This implies that the system of BSDEs with jumps (\ref{mainbsde}) turns into a decoupled one since the equations in (\ref{bsdeinter}) are not related each other.

Next for any $i=1,\dots,m$, the functions $u^i$, $g^i$ and $(t,x)\mapsto f^{(i)}(t,x,0,0,0)$ are of polynomial growth and finally $y\mapsto f^{(i)}(t,x,y,0,0)$ is Lipschitz uniformly $\wr$ $(t,x)$. Then for some $C$ and $\rho\geq 0$\be\label{estimfo}\E\Big [|g^i(X^{t,x}_T)|^p+(\int_0^T|f^{(i)}(r,X^{t,x}_r, (u^j(r,X^{t,x}_r))_{j=1,m},0,0)|^2dr)^{\frac{p}{2}}\Big ]\leq C(1+|x|^\rho).\ee

Let us now fix $i_0\in \mo$. Let $\cB^p$ be the space of processes $(Z,U)=(Z_s,U_s)_{s\leq T}$ such that:
\ms 

\no (a) $Z$ is $\cP$-measurable, $\r^{d}$-valued and $\E[(\int_0^T|Z_s|^2ds)^{\frac{p}{2}}]<\infty$ ;

\no (b) $U$ is $\bp$-measurable, $L^2(\l)$-valued and $\E[(\int_0^T\|U_s\|_{L^2(\l)}^2ds)^{\frac{p}{2}}]<\infty$. 
\ms

\no For $(\eta,\z)\in \cB^p$ let $\Phi(\eta,\z)=(\bar Z,\bar U)$ where $(\bar Y,\bar Z,\bar U)$ is the solution of the following BSDE: 
\be\label{edsrinter1}\left\{\ba{l}
\bar Y \in \cS^2(\r), \bar Z \in \cH^2(\r^{d}), \bar U\in \cH^2(L^2(\l));\\\\
\bar Y_s=g^{i_0}(X^{t,x}_T)+
 \int_s^Tf^{(i_0)}(r,X^{t,x}_r, (u^j(r,X^{t,x}_r))_{j=1,m},\eta_r, \z_r)dr -\int_s^T\bar Z_r dB_r-\int_s^T\int_{E} 
    \bar U_{r}(e)\tilde{\mu}(dr,de), \,\forall s\leq T.\ea \right.
\ee It implies that for any $s\le T$, 
\be\label{edsrinter2}
\bar Y_s=\E\Big [g^{i_0}(X^{t,x}_T)+
 \int_s^Tf^{({i_0})}(r,X^{t,x}_r, (u^j(r,X^{t,x}_r))_{j=1,m},\eta_r, \z_r)dr|\cF_s\Big ].
\ee
and then by Doob's martingale inequality and Jensen's one we deduce that
\be\label{ineqxx1}
\E\Big [\sup_{s\leq T}|\bar Y_s|^p\Big ]\leq C_p\E\Big [|g^{i_0}(X^{t,x}_T)|^p+T^{\frac{p}{2}}
 (\int_0^T|f^{({i_0})}(r,X^{t,x}_r, (u^j(r,X^{t,x}_r))_{j=1,m},\eta_r, \z_r)|^2dr)^{\frac{p}{2}}\Big ]
\ee
where $C_p$ is, along with this proof, a  constant independent of $T$ which may change from line to line. On the other hand, by the Burkholder-Davis-Gundy inequality and Doob's martingale one (see e.g. \cite{revuzyor}) we have   
$$\begin{array}{ll}
\E\Big [(\int_0^T|\bar Z_r|^2dr+
\int_0^T\|\bar U_r\|_{L^2(\l)}^2dr)^{\frac{p}{2}}\Big ]&\leq C_p\E\Big [\sup_{t\leq T}|\int_0^t\bar Z_r dB_r+\int_0^t\int_{E} 
    \bar U_{r}(e)\tilde{\mu}(dr,de)|^p\Big ]\\\\
  {}& \le C_p \E\Big [\{\sup_{s\leq T}|\bar Y_s| +\int_0^T|f^{({i_0})}(r,X^{t,x}_r, (u^j(r,X^{t,x}_r))_{j=1,m},\eta_r, \z_r)|dr\}^{p}\Big ]\ea
$$
and taking into account (\ref{ineqxx1}) and once more Jensen's inequality we deduce that 
\be\label{estimxxphi}\begin{array}{l}\dis{
\E\Big [(\int_0^T|\bar Z_r|^2dr+
\int_0^T\|\bar U_r\|_{L^2(\l)}^2dr)^{\frac{p}{2}}\Big ]\leq C_p\E\Big [|g^{i_0}(X^{t,x}_T)|^p+T^{\frac{p}{2}}
 (\int_0^T|f^{(i_0)}(r,X^{t,x}_r, (u^j(r,X^{t,x}_r))_{j=1,m},\eta_r, \z_r)|^2dr)^{\frac{p}{2}}\Big ].}\ea
 \ee
It means that $\Phi (\eta, \zeta)\in \cB^p$, for any $(\eta, \zeta)\in \cB^p$. On the other hand, let us set 
$(\bar Z^1,\bar U^1)=\Phi(\eta^1,\z^1)$. Then $(\bar Y-\bar Y^1, \bar Z-\bar Z^1, \bar U-\bar U^1)$ verify the following BSDE: for any $s\leq T$, 
\be \label{bsdedif}\begin{array}{ll}
\bar Y_s-
\bar Y^1_s=
\int_s^T\{f^{(i_0)}(r,X^{t,x}_r, (u^j(r,X^{t,x}_r))_{j=1,m},\eta_r, \z_r)-
f^{(i_0)}(r,X^{t,x}_r, (u^j(r,X^{t,x}_r))_{j=1,m},\eta_r^1, \z^1_r)\}dr \\\qquad\qquad \qquad\qquad  -\int_s^T(\bar Z_r-\bar Z^1_r) dB_r-\int_s^T\int_{E} 
   ( \bar U_{r}(e)-\bar U^1_r(e))\tilde{\mu}(dr,de).\end{array}
\ee
As $f^{(i_0)}$ is Lipschitz then, in the same way as previously in considering the  BSDE (\ref{bsdedif}), we obtain: 
$$\begin{array}{l}\dis{
\E\Big [(\int_0^T|\bar Z_r-\bar Z^1_r|^2dr+
\int_0^T\|\bar U_r-\bar U_r^1\|_{L^2(\l)}^2dr)^{\frac{p}{2}}\Big ]\leq C_p  T^{\frac{p}{2}}\E\Big [
 (\int_0^T(|\eta_r-\eta^1_r|^2+\|\z_r-\z^1_r)\|_{L^2(\l)}^2)dr))^{\frac{p}{2}}\Big ].}\ea
$$
Now let  $\delta>0$. In considering the previous BSDEs (\ref{edsrinter1})-(\ref{bsdedif}) for $t\in [T-\delta, T]$ we obtain, in a similar way as previously,  
$$\begin{array}{l}\dis{
\Big (\E\Big [(\int_{T-\d}^T(|\bar Z_r-\bar Z^1_r|^2+
\|\bar U_r-\bar U_r^1\|_{L^2(\l)}^2)dr)^{\frac{p}{2}}\Big ]\Big )^{\frac{1}{p}}\leq C_p  \sqrt{\d}\Big (\E\Big [
 (\int_{T-\d}^T(|\eta_r-\eta^1_r|^2+\|\z_r-\z^1_r)\|_{L^2(\l)}^2)  dr))^{\frac{p}{2}}\Big ]\Big) ^{\frac{1}{p}}.}\ea
$$
Take $\d=(4C_p^2)^{-1}$, we obtain that $\Phi$ is a contraction when we restrict time $s$ to the interval $[T-\d,T]$. Then it has a fixed point which is 
nothing else but $(Z^{i_0;t,x},U^{i_0;t,x})$ since the solution of the BSDE (\ref{bsdeinter}) is unique on $[T-\d,T]$. 
\\ Let us define now $\|(\eta,\zeta)\|_{\d, p}$ ($(\eta,\zeta)\in {\cal B}^p)$ by:
$$
\|(\eta,\zeta)\|_{\d, p}:=\{\E\Big [(\int_{T-\d}^T(|\eta_r|^2+
\|\z_r\|_{L^2(\l)}^2)dr)^{\frac{p}{2}}\Big ]\}^{\frac{1}{p}}.$$ Next let us consider the sequence of processes of ${\cal B}^p$ defined by:
$$
(z^0,\z^0)=(0,0) \mbox{ and for }n\geq 1,\,\,(z^n,\zeta^n)=\Phi(z^{n-1},\z^{n-1}).
$$
It implies that 
$$\ba{ll}
\|\Phi(Z^{i_0;t,x},U^{i_0;t,x})-\Phi(z^n,\zeta^n)\|_{p,\delta} &=
\|(Z^{i_0;t,x},U^{i_0;t,x})-\Phi(z^n,\zeta^n)\|_{p,\delta}\\ &\leq \frac{1}{2}
\|(Z^{i_0;t,x},U^{i_0;t,x})-(z^n,\zeta^n)\|_{p,\delta}
\end{array}
$$
and then 
$$\ba{ll}
\|(Z^{i_0;t,x},U^{i_0;t,x})-(z^n,\zeta^n)\|_{p,\delta}\le \frac{1}{2^{n}}\|(Z^{i_0;t,x},U^{i_0;t,x})\|_{p,\delta}.\end{array}
$$
But since $\Phi$ is a contraction then we can easily show that 
$$\forall n\ge 1,\,\,
\|(z^n,\zeta^n)\|_{p,\delta}=\|\Phi(z^{n-1},\zeta^{n-1})\|_{p,\delta}\le 2
\|(z^1,\zeta^1)\|_{p,\delta}.
$$
Thus for any $n\ge 1$ we have 
$$
\|(Z^{i_0;t,x},U^{i_0;t,x})\|_{p,\d}\leq (\frac{2^{n+1}}{2^{n}-1})
\|(z^1,\z^1)\|_{p,\d}\le 4\|(z^1,\z^1)\|_{p,\d}.
$$
Next in the same way as in (\ref{estimxxphi}) we have 
\be\label{estimxxphi(ii)}\begin{array}{l}\dis{
\|(z^1,\z^1)\|_{p,\d}\leq C_p(\E[|g^{i_0}(X^{t,x}_T)|^p+\delta^{\frac{p}{2}}
 (\int_{T-\d}^T|f^{(i_0)}(r,X^{t,x}_r, (u^j(r,X^{t,x}_r))_{j=1,m},0, 0)|^2dr)^{\frac{p}{2}}])^{\frac{1}{p}}}\ea
 \ee 
and then by (\ref{estimfo}) we deduce that, for some non-negative constants $C$ and $\rho$,
$$
\|(z^1,\z^1)\|_{p,\d}\leq C(1+|x|^\rho)
$$
which implies 
$$
\E\Big [(\int_{T-\d}^T\|U^{i_0;t,x}_r\|_{L^2(\l)}^2dr)^{\frac{p}{2}}\Big ]\leq C(1+|x|^\rho).
$$
Next on $[t,T-\d]$ we have 
$$\ba{l}
Y^{i_0;t,x}_s=u^{i_0}(T-\d,X^{t,x}_{T-\d})+
 \int_s^{T-\d}f^{(i_0)}(r,X^{t,x}_r, (u^{j}(r,X^{t,x}_r))_{j=1,m},Z_r^{i_0;t,x}, U_r^{i_0;t,x})dr \\\\\qq\qq\qq -\int_s^{T-\d}Z_r^{i_0;t,x} dB_r-\int_s^{T-\d}\int_{E} 
    U_{r}^{i_0;t,x}(e)\tilde{\mu}(dr,de), \,\,t\leq s\le T-\d.\ea$$
The same calculations as previously lead to 
$$
\E\Big [(\int_{T-2\d}^{T-\d}\|U^{i_0;t,x}_r\|_{L^2(\l)}^2dr)^{\frac{p}{2}}\Big ]\leq C(1+|x|^\rho).
$$
since $u^{i_0}$, like $g^{i_0}$, is of polynomial growth. Repeating now this procedure on $[T-3\d,T-2\d]$, etc., and by (\ref{uzavant}) we obtain 
$$
\E\Big [\Big \{\int_{0}^{T}\|U^{i_0;t,x}_r\|_{L^2(\l)}^2dr\Big \}^{\frac{p}{2}}\Big ]\leq C(1+|x|^\rho).
$$
Finally since $i_0\in \mo$ is arbitrary we then obtain the estimate (\ref{up}). 

\begin{rem} The result of Lemma \ref{uplem} holds for functions $f^{(i)}$, $i=1,...,m$, satisfying the properties (\ref{condfi})-(a),(b) only independently of the structure condition (\ref{deffi}). \qed
\end{rem}

\begin{propo} \lb{lipschitzlocu}For any $i=1,\dots,m$, $u^i$ belongs to $\u$.
\end{propo}
\noindent {\bf Proof}: Let $x$ and $x'$ be elements of $\r^k$. Let $(\vec{Y}^{t,x}, \zx,\ux)$ (resp. $(\vec{Y}^{t,x'},Z^{t,x'},U^{t,x'})$) be the solution of the BSDE (\ref{mainbsde}) associated with $(f(s,\x_s,y,z,\z), g(\x_T))$ (resp. $(f(s,X^{t,x'}_s,y,z,\z), g(X^{t,x'}_T))$). Applying It\^o's formula to 
$| \vec{Y}^{t,x}_s - \vec{Y}_{s}^{t,x'}|^{2}$
between $s$ and $T$ and taking expectation yields: $\forall s\in [t,T]$, 
\be\lb{eqdeltay}  \begin{array}{ll}
&\E[| \vec{Y}_s^{t,x} - \vec{Y}_{s}^{t,x'}|^{2}+\int_{s}^{T}  \big\{|\Delta Z_r|^2 + \|\Delta U_r\|_{L^2(\l )}^{2}  \big\} dr ]
   \; \\\\&\qq\qq\qq\qq =\E[ |g(X_{T}^{t,x}) - g(X_{T}^{t,x'})|^{2} 
+ 2\int_{s}^{T}\langle(\vec{Y}_r^{t,x} - \vec{Y}_r^{t,x'} ),\;\Delta f(r)
\rangle dr]\\
 \end{array}
\ee
where the four processes $\Delta f(r)$, $\Delta Z_r $, $\Delta U_r(e)$, $\D Y(r)$ and $\D X_r$  are defined as follows: $\forall r\in [t,T]$, 
 $$\begin{array}{l} \Delta f(r):=(\Delta f^{(i)}(r))_{i=1,m}= \big(f^{(i)}(r, X_r^{t,x},\vec{Y}_r^{t,x}, Z_r^{i;t,x},U_r^{i;t,x}) -f^{(i)}(r, X_r^{t,x'},\vec{Y}_r^{t,x'}, Z_r^{i;t,x'},U_r^{i;t,x'})\big)_{i=1,m},\,\,\\\\ \Delta Z_r = Z_{r}^{t,x} -Z_{r}^{t,x'},\, \Delta U_r = U_{r}^{t,x} -U_{r}^{t,x'}, 
\Delta Y(r)=\vec{Y}_r^{t,x}-\vec{Y}_r^{t,x'}=({Y}_r^{j;t,x} - {Y}_r^{j;t,x'})_{j=1,m}\mbox{ and }\D X_r:=X^{t,x}_r-X^{t,x'}_r. \end{array}$$
($\langle .,.\rangle $ is the usual scalar product on $\mathbb{R}^m$). As for any $i\in \mo$, $g^i$ belongs to $\u$ and by (\ref{estimx}) and finally using Cauchy-Schwarz inequality to obtain: 
$$ \mathbb{E}[|g(X_{T}^{t,x}) - g(X_{T}^{t,x'})|^2]\le  C (1+|x|^p+|x'|^p)^2|x-x'|^2.$$
Therefore, we only need to deal with the other term of the right-hand side of (\ref{eqdeltay}), i.e.,
$$ 2\mathbb{E}[\int_{s}^{T}\langle{\big(\vec{Y}_r^{t,x} - \vec{Y}_r^{t,x'}\big), \Delta f(r)\rangle} dr].$$\\
Taking into account the expression of $f^{(i)}$ given by (\ref{deffi}) we then split $\Delta f(r)$ in the following way: for $r\leq T$,
$$
\Delta f(r)=(\Delta f^{(i)}(r))_{i=1,m}=
\Delta_1(r)+\Delta_2(r)+\Delta_3(r)+
\Delta_4(r)=(
\Delta_1^i(r)+\Delta_2^i(r)+\Delta_3^i(r)+
\Delta_4^i(r))_{i=1,m}
$$
where for any $i=1,\dots,m$, 
\[\begin{array}{ll}
\Delta_1^i(r) \; & =  h^{(i)}(r, X^{t,x}_r, \vec{Y}^{t,x}_r, Z^{i;t,x}_r, \ie \g_i(r,\x_r,e) U^{i;t,x}_r(e)\l(de))-\\{}&
\qq\qq\qq h^{(i)}(r, X^{t,x'}_r, \vec{Y}^{t,x}_r, Z^{i;t,x}_r, \ie \g_i(r,\x_r,e) U^{i;t,x}_r(e)\l(de));\\\Delta_2^i(r) \; & =  h^{(i)}(r, X^{t,x'}_r, \vec{Y}^{t,x}_r, Z^{i;t,x}_r, \ie \g_i(r,\x_r,e) U^{i;t,x}_r(e)\l(de))-\\{}&
\qq\qq\qq h^{(i)}(r, X^{t,x'}_r, \vec{Y}^{t,x'}_r, Z^{i;t,x}_r, \ie \g_i(r,\x_r,e) U^{i;t,x}_r(e)\l(de));\\
\Delta_3^i(r) \; & =  h^{(i)}(r, X^{t,x'}_r, \vec{Y}^{t,x'}_r, Z^{i;t,x}_r, \ie \g_i(r,\x_r,e) U^{i;t,x}_r(e)\l(de))-\\{}&
\qq\qq\qq h^{(i)}(r, X^{t,x'}_r, \vec{Y}^{t,x'}_r, Z^{i;t,x'}_r, \ie \g_i(r,\x_r,e) U^{i;t,x}_r(e)\l(de));
\\
\Delta_4^i(r) \; & =  h^{(i)}(r, X^{t,x'}_r, \vec{Y}^{t,x'}_r, Z^{i;t,x'}_r, \ie \g_i(r,\x_r,e) U^{i;t,x}_r(e)\l(de))-\\{}&
\qq\qq\qq h^{(i)}(r, X^{t,x'}_r, \vec{Y}^{t,x'}_r, Z^{i;t,x'}_r, \ie \g_i(r,X^{t,x'}_r,e) U^{i;t,x'}_r(e)\l(de)).
\end{array}\] 
As $h^{(i)}$ verifies (\ref{polygrtf}) and then by estimate (\ref{estimx}) and Cauchy-Schwarz inequality we have: 
\be\lb{part1}\begin{array}{lll}\dis{
\mathbb{E}[2\int_{s}^{T} \langle{\Delta Y(r),\; \Delta_{1} (r)\rangle} dr ]}&\le  &\; \dis{
 \mathbb{E}[\frac{1}{\epsilon}\int_{s}^{T}|\Delta Y(r)|^{2}dr + C^2\epsilon \int_{s}^{T} (1+|X^{t,x}_r|^p+|X^{t,x'}_r|^p)^2
|X^{t,x}_r-X^{t,x'}_r|^2dr]}\\
 \\
 & \le &\;  \; \dis{ \mathbb{E} [\frac{1}{\epsilon}\int_{s}^{T}|\Delta Y(r)|^{2}dr ]+ C^2\epsilon (1+|x|^{2p}+|x'|^{2p})|x -x'|^{2}.}
 \end{array}
\ee
Besides since $h^{(i)}$ is Lipschitz w.r.t $(y,z,q)$  then 
\be\lb{part2} \mathbb{E}[2\int_{s}^{T} \langle{\Delta Y(r), \, \Delta_{2}(r)\rangle} dr] \le \; 2 C\mathbb{E}[ \int_{s}^{T} |\Delta Y(r)|^{2}dr ]\ee and \be\lb{part3}
\mathbb{E}[2\int_{s}^{T}\langle{\Delta Y(r),\;\Delta_{3}(r)\rangle} dr ] \le \; \mathbb{E}[\frac{1}{\epsilon}\int_{s}^{T}|\Delta Y(r)|^{2}dr + C^2\epsilon \int_{s}^{T} |\Delta Z_r|^{2} dr].\ee
It remains to obtain a control of the last term. But for any $s\in [t,T]$ we have, 
$$\displaystyle{
 \mathbb{E}\Big [ 2\int_{s}^{T} \langle{\Delta Y(r),\; \Delta_4 (r)\rangle}dr \Big ]\le
 2C \mathbb{E}\Big [ \int_{s}^{T}dr|\Delta Y(r)|\times | \int_E \big(U^{t,x}_r(e)\gamma(r,X^{t,x}_r,e) -U^{t,x'}_{r}(e)\gamma(r,X^{t,x'}_r,e)\big) \l(de)| \Big ].}
 $$
Next by splitting the crossing terms as follows: 
$$
U^{t,x}_s(e)\gamma(s,X^{t,x'}_s,e) -U^{t,x'}_{s}(e)\gamma(s,X^{t,x'}_s,e)
=\Delta U_s (e) \gamma(s,X^{t,x}_s,e) +U_s^{t,x'}(e)(\gamma(s,X^{t,x}_s,e)-
\gamma(s,X^{t,x'}_s,e))$$ and setting $\D \g_s (e):= (\gamma(s,X^{t,x}_s,e)-
\gamma(s,X^{t,x'}_s,e))$ we obtain:
\be\lb{part4} \begin{array}{ll}
\displaystyle{
 \mathbb{E}\Big [ 2\int_{s}^{T} \langle{\Delta Y (r),\;\Delta_4 (r) \rangle}dr\Big ]} &\;
\displaystyle{ \le 2C\mathbb{E}\Big [ \int_{s}^{T}|\Delta Y(r)|\times \big( \int_{E}  (|U_r^{t,x'}(e) \Delta\gamma_r(e)| +|\Delta U_r(e) \gamma(r,X^{t,x}_r,e)|)\l(de)\big) dr \Big ] }\\ 
 \\
 &\;\displaystyle{ \le  \frac{2}{\epsilon}\mathbb{E}\Big [ \int_s^{T} |\Delta Y(r) |^{2}dr  \Big ]+ C^{2}\epsilon \mathbb{E}\Big [\int_{s}^{T}\big(\int_{E}|U^{t,x'}_r(e) \Delta\gamma_r(e)| \l(de)\big)^{2}dr\Big ] }\\{}&\dis {
\qq \qq \qq \qq + \,C^{2}\epsilon \mathbb{E}\Big [\int_{s}^{T}\big(\int_{E}|\Delta U_r(e) \gamma(r,X^{t,x}_r,e)| \l(de)\big )^{2}dr \Big ].}
 \end{array} \ee By Cauchy-Schwarz inequality, (\ref{estimx}) and (\ref{condgamma})-(ii),  and the result of Lemma \ref{uplem} it holds:
\be\lb{part5} \ba{ll}
 &\mathbb{E}\left[\int_{s}^{T}\big(\int_{E}|U^{t,x'}_r(e) \Delta\gamma_r(e)| \l(de)\big)^{2}dr\right ]\\\\&\leq 
 \mathbb{E}\left[\int_{s}^{T}dr\big(\int_{E}|U^{t,x'}_r(e)|^2\l(de)\big)\big(\int_{E}|\Delta\gamma_r(e)|^{2} \l(de)\big)\right ]\\\\
 &\leq C\mathbb{E}\left[\{1+
 \sup_{r\in [t,T]}|X^{t,x}_r|^{2p}+|X^{t,x'}_r|^{2p})\times
 \sup_{r\in [t,T]}|X^{t,x}_r-X^{t,x'}_r|^2\}\times
  \int_{s}^{T}dr\big(\int_{E}|U^{t,x'}_r(e)|^2\l(de)\big)\right]\\\\
{}&\le C\sqrt{\mathbb{E}\left[
\{1+
 \sup_{r\in [t,T]}|X^{t,x}_r|^{4p}+|X^{t,x'}_r|^{4p})
\sup_{r\in [t,T]}|X^{t,x}_r-X^{t,x'}_r|^4\right ]} \times \sqrt{ \E\Big [\{\int_{s}^{T}dr\big(\int_{E}|U^{t,x'}_r(e)|^2\l(de))\}^2\Big ]}\\\\
{}&\le C|x-x'|^2(1+|x'|^p+|x|^p)
   \ea \ee for some exponent $p$. On the other hand using once more Cauchy-Schwarz inequality and (\ref{condgamma})-(i) we get
\be\lb{part6} \ba{ll}
 \mathbb{E}\left[\int_{s}^{T}\big(\int_{E}|\Delta U_r(e) \gamma(r,X^{t,x}_r,e)| \l(de)\big)^{2}dr\right ]&\leq 
 \mathbb{E}\left[\int_{s}^{T}dr\big(\int_{E}|\Delta U_r(e)|^2\l(de)\big)\big(\int_{E}|\gamma(r,X^{t,x}_r,e)|^{2} \l(de)\big)\right ]\\\\
 &\leq C\mathbb{E}\left[\int_{s}^{T}dr\big(\int_{E}|\Delta U_r(e)|^2\l(de)\big)\right].
   \ea\ee
Taking now into account inequalities (\ref{part1})-(\ref{part6}) we obtain: 
$$
\label{itoformula}
  \begin{array}{ll}
  \displaystyle{ | \vec{Y}_s^{t,x} - \vec{Y}_{s}^{t,x'}|^{2} } &\;\displaystyle{ +
  \mathbb{E}\Big [\int_{s}^{T}  \Big [|\Delta Z_r|^2 +  \|\Delta U_r\|_{L^2(\l )}^{2}  \Big ] dr \Big ] }
\\
 &\; \displaystyle{ = \mathbb{E}\Big [|g(X_{T}^{t,x}) - g(X_{T}^{t,x'})|^{2} \Big ]
+ 2  \mathbb{E}\Big [\int_{s}^{T} \langle{\big(\vec{Y}_r^{t,x} - \vec{Y}_r^{t,x'}\big),\; \Delta f(r)\rangle} dr \Big ]}\\
&\,\,\le \displaystyle{ |x-x'|^2(1+|x|^p+|x'|^p)(C+C^2\epsilon+C^3\epsilon )+(\frac{3}{\epsilon}+2C)\mathbb{E}\Big [\int_{s}^{T}|\Delta Y(r)|^{2}dr\Big ] +C^2\epsilon \int_{s}^{T} |\Delta Z_r|^{2} dr]}
 \\
&\qq\qq \displaystyle{ +C^3\epsilon \mathbb{E}\left[\int_{s}^{T}dr\big(\int_{E}|\Delta U_r(e)|^2\l(de)\big)\right]}.
 \end{array}
$$
Choosing now $\epsilon$ small enough we deduce the existence of a constant $C\ge 0$ such that for any $s\in [t,T]$, 
$$
\begin{array}{ll}
  \displaystyle{ \E[| \Delta Y(s)|^{2}] \leq  C|x-x'|^2(1+|x'|^{2p}+|x|^{2p})+C\mathbb{E}\Big [\int_{s}^{T}|\Delta Y(r)|^{2}dr\Big ]}
 \end{array}
$$and by Gronwall's inequality this implies that for any $s\in [t,T]$, $$\E[| \Delta Y(s)|^{2}]\leq C|x-x'|^2(1+|x'|^{2p}+|x|^{2p}).$$ 
Finally in taking $s=t$ and considering (\ref{repre}) we obtain the desired result. 
\begin{rem}\label{bu}
For any $\r$-valued function $v$ which belongs to $\cU$, the quantity $B_iv$ defined in (\ref{defbk}) is well posed since the functions $\beta$ and $(\g_i)_{i=1,m}$ verify (\ref{cdbeta}) and (\ref{condgamma}) respectively. \qed
\end{rem}

We are now going to express the process $U^{t,x}$ of the BSDE (\ref{mainbsde}) by means of the functions $(u^i)_{i=1,m}$. This relation between $u^i$ and $U^{i;t,x}$ is a second crucial point in this paper. Actually we have:
\begin{prop}\lb{represdesui} For any $i=1,...,m$, $(t,x)\in \esp$, 
\be\lb{expui}
U^{i;t,x}_s(e)=u^i(s,X^{t,x}_{s-}+\beta(s,X^{t,x}_{s-},e))-
u^i(s,X^{t,x}_{s-}),\,\,d \P\otimes ds\otimes d \l-ae \mb{ on }\Omega\times [t,T]\times E.
\ee
\end{prop}
\no {\bf Proof}: First note that since the measure $\l$ is not finite, then we cannot use the same technique as in \cite{Hammorlais-mai2015} because $U^{i;t,x}$ is only square integrable and not necessarily integrable wrt $d \P\otimes ds\otimes d \l$. Therefore we first begin by truncating the L\'evy measure. 
\paragraph*{Step 1: Truncation of the L\'evy measure}${}$

\noindent For any $k \ge 1$, let us first introduce a new Poisson random measure $\mu_k$ 
(obtained from the truncation of $\mu$) and its associated compensator $\nu_k$ as follows:
 $$ \displaystyle{\mu_k(ds,de)=1_{\{|e|\ge \frac{1}{k}\}}\mu(ds,de) \;\;\textrm{and} \;\;\nu_k(ds, de)= \l_k(de) ds :=1_{\{|e|\ge \frac{1}{k}\}}\l(de)ds}
$$
which means that, as usual, $\tilde{\mu}_k(ds,de) := (\mu_k -\nu_k)(ds, de)$, is the associated random martingale measure.
The main point to notice is that $\l_k(E)<\infty$ since $\l$ 
integrates $(1\wedge |e|^2)_{e\in E}$.  

Next, let us introduce the process $\xk$ solving the following standard SDE of jump-diffusion type:
  \be \label{eq:diffusion2}\begin{array}{l} \xk_{s}= x + \int_{t}^{s} b(r,\xk_r)dr +\int_{t}^{s}\sigma(r, \xk_r)dB_r\\ \\
 \qq\qq+\int_{t}^{s}\int_{E}\beta(r, \xk_{r-}, e) \tilde{\mu}_k(dr,de),\; s\in [t,T];\,\xk_s=x \mbox{ if }  s\leq t. \end{array}\ee

\no Note that thanks to the assumptions on $b$, $\sigma$ and $\beta$ the process $\xk$ exists and is unique. Morever it satisfies the same estimates as in (\ref{estimx}) since $\l_k$ is just a truncation at the origin of $\l$ which integrates $\ife$. 

On the other hand let us consider the following Markovian BSDE with jumps (similar as BSDE (\ref{mainbsde})):
 \be \label{auxil-bsde-jumps}
   \left\{\begin{array}{l} \E[\sup_{s\leq T}
|\yk_s|^2]+\int_0^T\{|\zk_r|^2+\int_E |\uk_r(e)|^2 \l_k(de)\}]<\infty \,\,;\\
\\\yk_s = g(\xk_T)+ \;\; \int_s^T f_{\mu_k}(r,\xk_r,\yk_r, \zk_r,\uk_r)dr\\ \\
\quad \qq\qq\qq\qq - \int_{s}^{T}  \{\zk_r  dB_r + \int_{E}{}\uk_r(e)\tilde{\mu}_k(dr,de)\}
,\;s\leq T,\\
\end{array}\right.\ee 
with, for any $(t,x)\in \esp$, $y\in \r^m$, $z=(z_i)_{i=1,m} \in \r^{m\times d}$ and $\z=(\z_i)_{i=1,m}\in L^2_{m}(E,\l_k)$, 
$$\ba{l}f_{\mu_k}(t,x,y,z,\z)=(f^{(i)}_{\mu_k}(t,x,y,z_i,\z_i))_{i=1,m}=
(h^{(i)}(t,x,y,z_i,\ie \g_i(t,x,e))\z_i(e)\l_k(de)))_{i=1,m}.\ea
$$

First let us emphasize that this latter BSDE is related to the filtration $(\cF ^k_s)_{s\le T}$ generated by the Brownian motion and the independant random measure $\mu_k$. However this point does not raise major issues since for any $s\leq T$, $\cF ^k_s\subset \cF_s$ and thanks to the relationship between $\mu$ and $\mu_k$. 

Next by the properties of the functions $b$, $\sigma$, $\beta$ and assumptions (H1), (H2), (\ref{condgamma}) on the functions $g$, $h$ and $\g$ respectively, and according to Proposition \ref{existencegene} (see also \cite{TangLi94} or \cite{BarlesBuckPardoux}), there exists a unique triple  $(\yk, \zk, \uk )$
 solving (\ref{auxil-bsde-jumps}). In addition, since the setting is Markovian, then by Proposition \ref{solvisco} there also exists a function $u^{k}$ from $\esp$ into $\r^m$ of $\pgc$ such that 
\be \label{new-FKformula} \forall s\in [t,T],\,\,\yk_s := u^k(s, \xk_s),\, \P-a.s. \ee
 Moreover as in Proposition \ref{lipschitzlocu}, there exist positive constants $C$ and $p$ which do not depend on $k$ such that:
 \be\lb{ukinu}
 \forall t,x,x',\,\,|u^k(t,x)-u^k(t,x')|\leq C(1+|x|^p+|x'|^p)|x-x'|.
 \ee
Finally as $\l_k$ is finite then we have the following relationship between the process 
$\uk=(^k U^{i;t,x})_{i=1,m}$ and the determinstic functions $u^k=(u^k_i)_{i=1,m}$ (see \cite{Hammorlais-mai2015}, Proposition 3.1, pp.6): $\forall i=1,\dots,m$, 
$$
^k U^{i;t,x}_s(e)=u^k_i(s,\xk_{s-}+\beta(s,\xk_{s-},e))-
u^k_i(s,\xk_{s-}),\,\,d \P\otimes ds\otimes d \l_k-ae \mb{ on }\Omega\times [t,T]\times E.
$$ This is mainly due to the fact that $\uk$   belongs to $L^1\cap L^2 (ds\otimes d\P\otimes d\l_k)$ since $\l_k(E)<\infty$ and then we can split the stochastic integral w.r.t $\tilde \mu_k$ in (\ref{auxil-bsde-jumps}). Therefore for all $i=1,...,m$, \be \lb{expuk}
^k U^{i;t,x}_s(e)1_{\{|e|\ge \frac{1}{k}\}}=(u^k_i(s,\xk_{s-}+\beta(s,\xk_{s-},e))-
u^k_i(s,\xk_{s-}))1_{\{|e|\ge \frac{1}{k}\}},\,\,d \P\otimes ds\otimes d \l-ae \mb{ on }\Omega\times [t,T]\times E.
\ee
\no{\bf Step 2: Convergence of the auxiliary processes}
\ms

Let us now prove the following convergence results
\be \label{cv1}
\E[\sup_{s\leq T}|X^{t,x}_s-\xk_s|^2]\rightarrow_k 0\ee 
and 
\be \label{cv2}\ba{l}
\E[\sup_{s\leq T}|Y^{t,x}_s-\yk_s|^2]+\int_0^T|Z^{t,x}_s-\zk_s|^2ds +\int_0^Tds\int_E \l(de)|U_s^{t,x}(e)-\uk_s(e)1_{\{|e|\geq \frac{1}{k}\}}|^2]\rightarrow_k 0.\ea\ee
where $X^{t,x} $ and  ($Y^{t,x},\; Z^{t,x},\; U^{t,x} $) are respectively solutions of the SDE (\ref{eq:diffusion}) and BSDE with jumps (\ref{mainbsde}).
\ms

First let us prove (\ref{cv1}) which is rather standard but we give it for completeness. For any $s\in [0,T]$ we have: 
\[  \ba{ll}  \displaystyle{ X_s^{t,x}-\xk_s = \;} & \; 
 \displaystyle{ \int_{0}^{s} \big( b(r,X^{t,x}_r)-b(r,\xk_r) \big)dr +\int_{0}^{s}\big( \sigma(r,X^{t,x}_r)-\sigma(r,\xk_r)\big) dB_r }\\
  \; & \displaystyle{ \qq\qq + \;\int_{0}^{s} \int_{E} \big(\beta(r,X^{t,x}_{r-},e) -\beta(r,\xk_{r-},e)
 \mathbf{1}_{\{|e| \ge \frac{1}{k}\}}\big)\tilde{\mu}(de, dr)}.\\
\ea \]
Next let $\et \in [0,T]$. Since $|a+b+c|^2 \le 3(|a|^2+|b|^2+|c|^2)$ for any real constants $a,b$ and $c$ and by the Cauchy-Schwarz and Burkholder-Davis-Gundy inequalities we have: 
\[
 \ba{l}
 \displaystyle{ \mathbb{E}\Big\{\sup_{0\le s \le \et} |X_s^{t,x}-\xk_s|^2 \Big\}} \\  
 \\
 \quad  \le \; \displaystyle{ 3 \mathbb{E}\Big\{\sup_{0 \le s \le \et}  | \int_{0}^{s} (b(r,X^{t,x}_r)-b(r,\xk_r)) dr| ^2
+  \sup_{0 \le s \le \et}| \int_{0}^{s} (\sigma(r,X^{t,x}_r)-\sigma(r,\xk_r))dB_r |^2 } \\
  \displaystyle{\quad \;   \;\quad  \quad +  \sup_{0 \le s \le \et}|\int_{0}^{s} \int_{E} 
  \big( \beta(r,X^{t,x}_{r-},e) -\beta(r,\xk_{r-},e)\mathbf{1}_{\{|e| \ge \frac{1}{k}\}}\big)
 \tilde{\mu}(de, dr)) |^2\Big\} } \\
%  %\\
 \quad  \le \; \displaystyle{ C 
 \mathbb{E}\Big\{ \int_{0}^\et \sup_{0 \le \t \le r} \{ |b(\t,X^{t,x}_\t)-b(\t,\xk_\t) | ^2 +|\sigma(\t,X^{t,x}_\t)-\sigma(\t,\xk_\t)|^{2}\}dr\}\Big\} }\\
   \quad \quad  \displaystyle{ +\;
   C \mathbb{E}\Big\{\int_{0}^\et \int_E  \sup_{0 \le \t \le r} \big| \beta(\t,X^{t,x}_{\t-},e) -\beta(\t,\xk_{\t-},e)\big|^2 \l_k(de) dr +
   \int_{0}^\et\int_E  \sup_{0 \le \t \le r} | \beta(\t,X^{t,x}_{\t-},e) |^2 \mathbf{1}_{\{|e| < \frac{1}{k}\}}\l(de)dr\Big\}}.\ea 
\] But $b$ and $\sigma$ are Lipschitz $\wr$ $x$ and $\b$ verifies (\ref{cdbeta})-(ii), then we have: $\forall r\in [0,T]$, 
 $$\displaystyle{ \sup_{0 \le \t \le r} \{ |b(\t,X^{t,x}_\t)-b(\t,\xk_\t) | ^2 +|\sigma(\t,X^{t,x}_\t)-\sigma(\t,\xk_\t)|^{2}\} \le 
 C \sup_{0\le \t \le r} |X^{t,x}_{\t} -\xk_{\t}|^2}$$
 and 
 $$\begin{array}{ll}\int_E \sup_{0 \le \t \le r} \big| \beta(\t,X^{t,x}_{\t-},e) -\beta(\t,\xk_{\t-},e)\big|^2 \l_k(de) 
 &\le C \sup_{0\le \t \le r} |X^{t,x}_{\t} -\xk_{\t}|^2 \int_{E} (1 \wedge|e|)^2\l_k(de) \\\\
{}&\le C \sup_{0\le \t \le r} |X^{t,x}_{\t} -\xk_{\t}|^2  
 \ea$$
for some constant $C$ since $\l_k((1\wedge |e|^2)_{e\in E})$ is smaller than $\l((1\wedge |e|^2)_{e\in E})$ and this quantity is  finite. Plug now those two last inequalities in the previous one to obtain: $\forall \et \in [0,T]$,
$$
 \ba{l}\displaystyle{
\mathbb{E}\Big\{\sup_{0\le s \le \et} |X_s^{t,x}-\xk_s|^2 \Big\}\le  C\E\Big\{\int_0^\eta
  \sup_{0\le \t \le r} 
  |X^{t,x}_{\t} -\xk_{\t}|^2 dr +
 C\int_{\{|e| < \frac{1}{k}\}}(1\wedge |e|^2)
 \l(de)\Big\}.}
 \end{array}$$ Finally by Gronwall's Lemma we obtain the desired result since $\int_{\{|e| < \frac{1}{k}\}}(1\wedge |e|^2)
 \l(de)\rwk 0$.
 
 We now focus on (\ref{cv2}). First note that we can apply It\^o's formula, even if the BSDEs are related to filtrations and Poisson random measures which are not the same, since:
  
 (i) $\cF^k_s\subset \cF_s$, $\forall s\leq T$ ;
 
 (ii) for any $s\leq T$, $\int_0^s\int_E {^k\!U}^{i;t,x}_r(e)\tilde \mu_k(dr,de)=
 \int_0^s\int_E{^k\!U}^{i;t,x}_r(e){1}_{\{|e|\geq \frac{1}{k}\}}\tilde \mu(dr,de)$ and then the first $(\cF^k_s)_{s\leq T}$-martingale is also an 
 $(\cF_s)_{s\leq T}$-martingale.
\ms

\no Therefore we have: $\forall s\in [0,T]$,
\be\lb{itoyyk}\ba{l}
\E[|\vec{Y}_s^{t,x}-\yk_s|^2+\int_0^T\{|\zx_s-\zk_s|^2+\int_E|\ux_s(e)-\uk_s(e)1_{\{|e|\ge \frac{1}{k}\}}|^2\l(de)\}ds]=\E[|g(X_T^{t,x})-g(\xk_T)|^2]\\\\\qq\qq\qq +2\E[\int_s^T
(\vec{Y}_r^{t,x}-\yk_r)\times( f(r,\x_r,\vec{Y}_r^{t,x},\zx_r,\ux_r)-
f_k(r,\xk_r,\yk_r,\zk_r,\uk_r))dr].
\ea
\ee
First note that by (\ref{cv1}) and since $g$ belongs to $\u$ and $\xk$ verifies estimates (\ref{estimx}) then it holds:\be \lb{cvtermg}\E[|g(X_T^{t,x})-g(\xk_T)|^2]\rwk 0.\ee Next let us set: 
$$
(f(r,\x_r,\vec{Y}_r^{t,x},\zx_r,\ux_r)-
f_k(r,\xk_r,\yk_r,\zk_r,\uk_r))=A(r)+B(r)+C(r)+D(r)
$$
where, taking into account the expression of $f$ through $h$ (see (\ref{deffi})), for any $r\in [0,T]$: 
$$\ba{l}
\mb{(i)}\q A(r)=(h^{(i)}(r,\x_r,\vec{Y}_r^{t,x},Z^{i;t,x}_r,\ie \g_i(r,\x_r,e)U^{i;t,x}_r(e)\l(de))-\\\qq \qq \qq \qq \qq \q h^{(i)}(r,\xk_r,\vec{Y}_r^{t,x},Z^{i;t,x}_r,\ie \g_i(r,\x_r,e)U^{i;t,x}_r(e)\l(de)))_{i=1,m}\,\,;\\\\
\mb{(ii)}\qq B(r)=(h^{(i)}(r,\xk_r,\vec{Y}_r^{t,x},Z^{i;t,x}_r,\ie \g_i(r,\x_r,e)U^{i;t,x}_r(e)\l(de))-\\\qq \qq \qq \qq \qq \qq h^{(i)}(r,\xk_r,\yk_r,Z^{i;t,x}_r,\ie \g_i(r,\x_r,e)U^{i;t,x}_r(e)\l(de)))_{i=1,m}\,\,;\\\\
\mb{(iii)} \q
C(r)=(h^{(i)}(r,\xk_r,\yk_r,Z^{i;t,x}_r,\ie \g_i(r,\x_r,e)U^{i;t,x}_r(e)\l(de))-\\\qq \qq \qq \qq \qq \qq h^{(i)}(r,\xk_r,\yk_r,^k\!\! Z^{i;t,x}_r,\ie \g_i(r,\x_r,e)U^{i;t,x}_r(e)\l(de)))_{i=1,m}\,\,;\\\\
\mb{(iv)}\q D(r)=(h^{(i)}(r,\xk_r,\yk_r,^k\!\! Z^{i;t,x}_r,\ie \g_i(r,\x_r,e)U^{i;t,x}_r(e)\l(de))-\\\qq\qq\qq\qq\qq
h^{(i)}(r,\xk_r,\yk_r,^k\!\! Z^{i;t,x}_r,\ie \g_i(r,\xk_r,e)^k\!U^{i;t,x}_r(e)\l_k(de)))_{i=1,m}.\ea$$But by (\ref{cdlipschitzfyzu}) and 
(\ref{polygrtf}), we have: $\forall r\in [0,T]$,  \be \lb{splitf}\ba{l}
|A(r)|\leq C(1+|\x_r|^p+|\xk_r|^p)|\x_r-\xk_r|,\,\,|B(r)|\leq C|\vec{Y}_r^{t,x}-\yk_r| \mb{ and }
|C(r)|\le |\zx_r-\zk_r|\ea\ee
 where $C$ is a constant. Finally let us deal with $D(r)$ which is more involved. First note that $D(r)=(D_i(r))_{i=1,m}$ where 
$$\ba{l}
D_i(r)=
h^{(i)}(r,\xk_r,\yk_r,\zki_r,\ie\g_i(r,\x_r,e)U^{i;t,x}_r(e)\l(de))-\\\qq\qq\qq\qq 
h^{(i)}(r,\xk_r,\yk_r,\zki_r,\ie \g_i(r,\xk_r,e){^k\! U^{i;t,x}_r(e)}\l_k(de)).\ea
$$But as $h^{(i)}$ is Lipschitz w.r.t to the last component $q$ then 
$$\ba{ll}
|D_i(r)|^2&\le C\{\ie |\g_i(r,\x_r,e)U^{i;t,x}_r(e)- \g_i(r,\xk_r,e){^k\! U^{i;t,x}_r(e)}\itk|\l(de)\}^2\\\\
{}&\le C\Big\{\{\ie |\g_i(r,\x_r,e)-\g_i(r,\xk_r,e)||U^{i;t,x}_r(e)|\l(de)\}^2\\&\qq\qq\qq\qq+\{\ie|\g_i(r,\xk_r,e)||
U^{i;t,x}_r(e)- {^k\! U^{i;t,x}_r(e)}\itk|\l(de)\}^2\Big \}\\\\&\le C\{(1+|\x_r|^p+|\xk_r|^p)|\x_r-\xk_r|\ie (1\wedge |e|)|U^{i;t,x}_r(e)|\l(de)\}^2\\&\qq\qq\qq\qq+C
\ie |
U^{i;t,x}_r(e)- {^k\! U^{i;t,x}_r(e)}\itk|^2\l(de).
\ea
$$The last inequality follows from the properties (\ref{condgamma})-(i), (ii) satisfied by $\g_i$ and Cauchy-Schwarz inequality. Next going back to
(\ref{itoyyk}) and arguing as in the bulk of the proof of Proposition \ref{lipschitzlocu} we deduce the existence of a constant $C\ge 0$ independant of $k$ such that:
\be\lb{itoyyk2}\ba{l}
\E[|\vec{Y}_s^{t,x}-\yk_s|^2+\int_s^T\{|\zx_s-\zk_s|^2+\int_E|\ux_s(e)-\uk_s(e)1_{\{|e|\ge \frac{1}{k}\}}|^2\l(de)\}ds]\\\\\le C\E[|g(X_T^{t,x})-g(\xk_T)|^2]+C\E[\int_s^T|\vec{Y}_r^{t,x}-\yk_r|^2dr]
+C\E[\int_0^T
(1+|\x_r|^p+|\xk_r|^p)^2|\x_r-\xk_r|^2dr]\\\\\qq +C\E[\int_0^Tdr
\{(1+|\x_r|^p+|\xk_r|^p)|\x_r-\xk_r|\ie (1\wedge |e|)|U^{i;t,x}_r(e)|\l(de)\}^2], \forall s\in [0,T].
\ea
\ee
But $$\ba{l}\E[|g(X_T^{t,x})-g(\xk_T)|^2]+\E[\int_0^T
(1+|\x_r|^p+|\xk_r|^p)^2|\x_r-\xk_r|^2dr]\rwk 0\ea$$ and $$\ba{l}\E[\int_0^Tdr
\{(1+|\x_r|^p+|\xk_r|^p)|\x_r-\xk_r|\ie (1\wedge |e|)|U^{i;t,x}_r(e)|\l(de)\}^2]\rwk 0.\ea$$ 
Let us focus indeed on the first convergence. Obviously the first term converges to 0 because $g$ belongs to $\u$ and $X^{t,x}$, $\xk$ verify estimates (\ref{estimx}) uniformly and by (\ref{cv1}). For the second term we have: 
$$\begin{array}{ll}
&\E[\int_0^T
(1+|\x_r|^p+|\xk_r|^p)^2|\x_r-\xk_r|^2dr]\\{}&\qquad \qquad \le \E[\sup_{r\leq T}|\x_r-\xk_r|\int_0^T
(1+|\x_r|^p+|\xk_r|^p)^2|\x_r-\xk_r|dr]\\
{}&\qquad \qquad \le  \{\E[\sup_{r\leq T}|\x_r-\xk_r|^2]\}^{\frac{1}{2}}\{\E[(\int_0^T
(1+|\x_r|^p+|\xk_r|^p)^2|\x_r-\xk_r|dr)^2]\}^{\frac{1}{2}}.
\ea
$$
But the first factor in the right-hand side of this inequality goes to 0 when $k\rw \infty$ due to (\ref{cv1}) and the second factor is uniformly bounded by the uniform estimates (\ref{estimx}) of $X^{t,x}$ and $\xk$. \\
For the second convergence, it is a consequence of (\ref{cv1}), the fact that $\xk$ verifies estimates (\ref{estimx}) uniformly, the Cauchy-Schwarz inequality (used twice) and finally (\ref{up}). Then by Gronwall's Lemma we deduce first that for any $s\leq T$, $$
\E[|\vec{Y}_s^{t,x}-\yk_s|^2]\rwk 0
$$
and in taking $s=t$ we obtain $u^k(t,x)\rwk u(t,x)$. As $\tx$ is arbitrary then $u^k \rwk u$ pointwisely. Next going back to (\ref{itoyyk2}) take the limit w.r.t $k$ and using the uniform polynomial growth of $u^k$ and the Lebesgue dominated convergence theorem as well, to obtain:
\be\lb{limuk}\E[\int_t^T\int_E|\ux_s(e)-\uk_s(e)1_{\{|e|\ge \frac{1}{k}\}}|^2\l(de)\}ds]\rwk 0.\ee
\no{\bf Step 3: Conclusion}
\ms

\no First note that by (\ref{ukinu}) and the pointwise convergence of $(u^k)_k$ to $u$, if $(x_k)_k$ is a sequence of $\r^k$ which converges to $x$ then 
$(u^k(t,x_k))_k$ converges to $u(t,x)$. Now let us consider a subsequence which we still denote by $\{k\}$ such that $\sup_{s\leq T}|\xk_s-\x_s|^2\rwk 0$, $\P-a.s.$ (and then 
$|\xk_{s-}-\x_{s-}|\rw_k 0$ since $|\xk_{s-}-\x_{s-}|\leq \sup_{s\leq T}|\xk_s-\x_s|^2$). By (\ref{cv1}), this subsequence exists. As the mapping $x\mapsto \beta(t,x,e)$ is Lipschitz then the sequence 
\be\lb{convxx}\begin{array}{l}(
^k U^{i;t,x}_s(e)1_{\{|e|\ge \frac{1}{k}\}})_k=
((u^k_i(s,\xk_{s-}+\beta(s,\xk_{s-},e))-
u^k_i(s,\xk_{s-}))1_{\{|e|\ge \frac{1}{k}\}})_{k\geq 1}\rw_k\\\\
\qq\qq\qq (u_i(s,\x_{s-}+\beta(s,\x_{s-},e))-
u_i(s,\x_{s-})), d \P\otimes ds\otimes d \l-ae \mb{ on }\Omega\times [t,T]\times E.\ea\ee for any $i=1,...,m$. Finally from (\ref{limuk}) we deduce that 
$$U^{t,x}_s(e)=(u(s,\x_{s-}+\beta(s,\x_{s-},e))-
u(s,\x_{s-})),\,\,d \P\otimes ds\otimes d \l-ae \mb{ on }\Omega\times [t,T]\times E
$$
which is the desired result. 
\begin{rem}\lb{remimpo}In order to prove the final step we do not need to use the property (\ref{ukinu}) satisfied by $u^k_i$. Instead, we only need that for 
any sequence $(x_k)_k$ which converges to $x$, the sequence $(u^k_i(t,x_k)-u^k_i(t,x))_k$ converges to $0$ and $(u^k_i(t,x))_{k\geq 1}$ converges to 
$u_i(t,x)$ pointwisely. This point plays an important role in the proof of uniqueness of Theorem \ref{mainresult}. \qed
\end{rem}

\section{The main result}
We are now ready to give the main result of this paper. Before doing so we recall the notion of viscosity solution we deal with. This definition has been more or less introduced in \cite{Hammorlais-mai2015}. 

For $\phi \in \cC^{1,2}(\esp)$, let us denote by 
$\cL^X\phi(t,x)$ the differential-integral generator associated with the jump-diffusion process introduced in (\ref{eq:diffusion}) and which is given by:  $\forall (t,x)\in \esp$, $$\ba{l}
\cL^X\phi(t,x):=
b(t,x)^\top D_x\phi(t,x)+\\\qq\qq\frac{1}{2}\mathrm{Tr}\big(\sigma\sigma^\top (t,x) D^2_{xx}\phi(t,x)\big)+\ie\{ \phi(t,x+\beta(t,x,e))-
\phi(t,x)-\beta(t,x,e)^\top D_x\phi(t,x)\}\l(de).\ea$$
\begin{defn}\label{od}
A family of deterministic functions $u =(u^i)_{i=1,m}$, such that, for any $i\in \mo$ $u^i: (t,x) \in \esp \mapsto u^i(t,x)\in \r$ belongs to the class $\mathcal{U}$, is said to be a viscosity sub-solution (resp. super-solution)
 of the IPDE (\ref{ipde-intro}) if: $\forall i\in \mo$, \\
(i) $\forall x\in \r^k$, $u^i(T, x) \le g^i(x) $ (resp. $u^i(T, x) \ge g^i(x) $) ;\\
 (ii)  For any  $(t,x)\in (0,T)\times \r^k$ and any function $\phi$ of class $\cC^{1,2}(\esp)$ such that $(t,x)$ is a global maximum (resp. minimum) point of 
$ u^i -\phi$ and $(u^i-\phi)(t,x)=0$, one has
$$ -\partial_t \phi(t,x)
-\mathcal{L}^X \phi(t,x) -h^{(i)}(t,x, (u^j(t,x))_{j=1,m}, \sigma^\top(t,x)D_x\phi(t,x), B_iu^{i}(t,x)) \le 0, 
$$
(resp. 
$$-\partial_t \phi(t,x)
-\mathcal{{L}}^X \phi(t,x) -h^{(i)}(t,x, (u^j(t,x))_{j=1,m}, \sigma^\top(t,x)D_x\phi(t,x), B_iu^{i}(t,x)) \ge 0).
$$

The family $u =(u^i)_{i=1,m}$ is a viscosity solution of (\ref{ipde-intro}) if it is both a viscosity sub-solution and viscosity super-solution. 
\end{defn}

Let us mention here the main difference with the classical definition of viscosity solution of (\ref{ipde-intro}) by Barles et al. \cite{BarlesBuckPardoux} (see Definition \ref{bbpdef} in Appendix). In our definition we keep $B_i u^i(t,x)$ which is defined since $u^i\in \u$ while in \cite{BarlesBuckPardoux} it is replaced with $B_i\phi(t,x)$ where  $\phi$ is the test function. This is one of the main reasons for which in \cite{BarlesBuckPardoux}, the authors have required monotonicity conditions (a)-(b) related to the functions $(\g_i)_\im$ and $(h^{(i)})_\im$. On the other hand note that when,  for any $i=1,...,m$, $h^{(i)}$ does not depend on $B_i u^i(t,x)$ those definitions coincide. 
\ms
 
We are now ready to state the main result of this paper.
\begin{thm} \label{mainresult}
Assume that Assumptions (H1)-(H2) are fulfilled. Then the $m$-tuple of functions $(u^i)_{i=1,m}$ defined in (\ref{repre}) is the unique viscosity solution of system (\ref{ipde-intro}) according to Definition \ref{od}. 
\end{thm}
\proof \underline{Step 1}: Existence 

Let us consider the following multi-dimensional BSDE:
\be \left\{ \label{mainbsde2}
    \begin{array}{l}  
\vec{\underbar Y}^{t,x}:=(\underbar Y^{i;t,x})_{i=1,m}\in \cS^2(\r^m),\,
{\underbar Z}^{t,x}:=(\underbar Z^{i;t,x})_{i=1,m}\in \cH^2(\r^{m\times d}),\, \underbar U^{t,x}:=(\underbar U^{i;t,x})_{i=1,m}\in \cH^2(L_m^2(\l));\\\\ \forall i\in \mo,\,\,
\underbar Y_T^{i;t,x}= g^{i}(X^{t,x}_T) \mbox{ and }\forall s\le T, \\\\
\q d\underbar Y_s^{i;t,x} = -h^{(i)}(s,X^{t,x}_s, \vec{\underbar Y}^{t,x}_s,\underbar Z_s^{i;t,x}, \int_E
\g_i(s,X^{t,x}_s,e)\{u^i(s,X^{t,x}_{s-}+\beta (s,X^{t,x}_{s-},e))-u^i(s,X^{t,x}_{s-})\}\l(de))ds\\\\ \qq\qq\qq\qq +\underbar Z_s^{i;t,x} dB_s +\int_{E} 
    \underbar U_{s}^{i;t,x}(e)\tilde{\mu}(ds,de).\end{array}\right.  \ee
Since for any $i=1,...,m$, $u^i$  belongs to $\u$, $\beta(t,x,e)$ and $\g_i(t,x,e)$ verify respectively (\ref{cdbeta}) and (\ref{condgamma})  and finally by Assumption (H2) we have:

(i) the mapping $(y,z)\mapsto h^{(i)}(s,\x_s,y,z,\int_E
\g_i(s,X^{t,x}_s,e)\{u^i(s,X^{t,x}_{s-}+\beta (s,X^{t,x}_{s-},e))-u^i(s,X^{t,x}_{s-})\}\l(de))$ is uniformly Lipschitz ;

(ii) the process $(h^{(i)}(s,\x_s,0,0,\int_E
\g_i(s,X^{t,x}_s,e)\{u^i(s,X^{t,x}_{s-}+\beta (s,X^{t,x}_{s-},e))-u^i(s,X^{t,x}_{s-})\}\l(de)))_{s\leq T}$ is $ds\otimes d\P$-square integrable. 
\ms

\noindent It follows that the solution of this backward equation (\ref{mainbsde2}) exists and is unique by Proposition \ref{existencegene} (see Remark \ref{remexistuniq}). Moreover, as the process $X^{t,x}$ is RCLL then the set of its discontinuous points on $[0,T]$ is at most countable. Therefore $\P-a.s., \mb{ for any }s\leq T$, it holds 
$$\begin{array}{l}
\int_s^Th^{(i)}(r,X^{t,x}_r, \vec{\underbar Y}^{t,x}_r,\underbar Z_r^{i;t,x}, \int_E
\g_i(r,X^{t,x}_r,e)\{u^i(r,X^{t,x}_{r-}+\beta (r,X^{t,x}_{r-},e))-u^i(r,X^{t,x}_{r-})\}\l(de))dr=\\\qq\qq\int_s^T
h^{(i)}(r,X^{t,x}_r, \vec{\underbar Y}^{t,x}_r,\underbar Z_r^{i;t,x}, \int_E
\g_i(r,X^{t,x}_r,e)\{u^i(r,X^{t,x}_{r}+\beta (r,X^{t,x}_{r},e))-u^i(r,X^{t,x}_{r})\}\l(de))dr.\ea $$
Next as for any $i=1,...,m$, $u^i$ belongs to $\u$, then by Proposition \ref{solvisco}, there exists a family of deterministic continuous functions of polynomial growth $(\underbar u^i)_{i=1,m}$ such that for any $(t,x)\in \esp$,
$$
\forall s\in [t,T],\,\, \underbar Y^{i;t,x}_s=
\underbar u^i(s,X^{t,x}_s).
$$
Finally, again by Proposition \ref{solvisco}, the family $(\underbar u^i)_{i=1,m}$ is a viscosity solution of the following system:
\be \left\{ \label{secondorder-pde-inter} \begin{array}{l} 
-\partial_t \underbar u^i(t,x)- b(t,x)^\top D_x\underbar u^i(t,x)-\frac{1}{2}\mathrm{Tr}\big(\sigma\sigma^\top (t,x) D^2_{xx}\underbar u^i (t,x)\big)-K\underbar u^i(t,x)\\\q\qq -h^{(i)}(t,x, (\underbar u^j(t,x))_{j=1,m}, (\sigma^\top D_x\underbar u^i)(t,x), B_iu^i(t,x))=0,\,\, (t,x)\in \esp ;\\
\underbar u^i(T,x) = g^i(x)
\end{array}  \right.
\ee
Note that in this system (\ref{secondorder-pde-inter}), the last component of $h^{(i)}$ is $B_iu^i(t,x)$ and not $B_i\underbar u^i(t,x)$. Next and once more, let us consider the system of BSDEs by which the family $(u^i)_{i=1,m}$ is defined through the Feynman Kac's formula (\ref{repre}): 
\be \left\{ \label{mainbsde21}
    \begin{array}{l}  
\vec{ Y}^{t,x}:=( Y^{i;t,x})_{i=1,m}\in \cS^2(\r^m),\,
{ Z}^{t,x}:=
(Z^{i;t,x})_{i=1,m} \in \cH^2(\r^{m\times d}),\,  U^{t,x}:=( U^{i;t,x})_{i=1,m}\in \cH^2(L_m^2(\l));\\\\ \forall i\in \mo,\,\,
 Y_T^{i;t,x}= g^{i}(X^{t,x}_T) \mbox{ and }\forall s\le T,\\\\
\q dY_s^{i;t,x} = -h^{(i)}(s,X^{t,x}_s, \vec{Y}^{t,x}_s,Z_s^{i;t,x}, \int_E
\g_i(s,X^{t,x}_s,e) U_{s}^{i;t,x}(e)\l(de))ds +Z_s^{i;t,x} dB_s +\int_{E} 
    U_{s}^{i;t,x}(e)\tilde{\mu}(ds,de).\end{array}\right.\ee
But by Proposition \ref{represdesui} we now 
that for any $i=1,m$, 
$$
U_{s}^{i;t,x}(e)=
u^i(s,X^{t,x}_{s-}+\beta (s,X^{t,x}_{s-},e))-u^i(s,X^{t,x}_{s-}),\,\,ds\otimes d\mathbb{P}\otimes d\l \mbox{ on }[t,T]\times \Omega \times E.
$$
Plug now this relation in the first term of the right-hand side of the second equality of (\ref{mainbsde21}), one obtains, by uniqueness of the solution of the BSDE (\ref{mainbsde2}), that for any $s\in [t,T] \mbox{ and }i\in \mo,\,\,
\underbar Y^{i;t,x}_s=Y^{i;t,x}_s$. Thus for any $i\in \mo$, $u^i=\underbar u^i$. Henceforth, the family $(u^i)_{i=1,m}$ is a viscosity solution of (\ref{ipde-intro}) in the sense of Definition \ref{od}. \\
\ms

\no \underline{Step 2}: Uniqueness
\ms

We now show uniqueness of the solution in the class $\u$. So let $(\bar u^i)_{i=1,m}$ be another family of $\u$ which is solution of the system (\ref{ipde-intro}) in the sense of Definition \ref{od} and let us consider the following system of BSDEs:
\be \left\{ \label{mainbsde3}
    \begin{array}{l}  
\vec{\bar Y}^{t,x}:=(\bar Y^{i;t,x})_{i=1,m}\in \cS^2(\r^m),\,
{\bar Z}^{t,x}:=(
\bar Z^{i;t,x})_{i=1,m} \in \cH^2(\r^{m\times d}),\, \bar U^{t,x}:=(\bar U^{i;t,x})_{i=1,m}\in \cH^2(L_m^2(\l));\\\\ \forall i\in \mo,\,\,
\bar Y_T^{i;t,x}= g^{i}(X^{t,x}_T) \mbox{ and }\forall s\le T,\\\\
\qq d\bar Y_s^{i;t,x} = -h^{(i)}(s,X^{t,x}_s, \vec{\bar Y}^{t,x}_s,\bar Z_s^{i;t,x}, \int_E
\g_i(s,X^{t,x}_s,e)\{\bar u^i(s,X^{t,x}_{s-}+\beta (s,X^{t,x}_{s-},e))-\bar u^i(s,X^{t,x}_{s-})\}\l(de))ds\\\\ \qq\qq\qq +\bar Z_s^{i;t,x} dB_s +\int_{E} 
    \bar U_{s}^{i;t,x}(e)\tilde{\mu}(ds,de).\end{array}\right.  \ee
As for the BSDE (\ref{mainbsde2}), the solution of the BSDE (\ref{mainbsde3}) exists and is unique since $(\bar u^i)_{i=1,m}$ belong to $\u$. Moreover there exists a family of deterministic continuous functions $(v^i)_{i=1,m}$ of class $\Pi_g$ such that 
$$
\forall s\in [t,T],\,\, \bar Y^{i;t,x}_s=v^i(s,X^{t,x}_s).
$$
Additionally, by Proposition \ref{solvisco}, $(v^i)_{i=1,m}$ is the unique solution in the subclass $\pgc$ of continuous functions of the following system: $\forall i=1,\dots,m$,  

\be \left\{ \label{secondorder-pde-inter2} \begin{array}{l} 
-\partial_t v^i(t,x)- b(t,x)^\top D_xv^i(t,x)-\frac{1}{2}\mathrm{Tr}\big(\sigma\sigma^\top (t,x) D^2_{xx}v^i(t,x)\big)-Kv^i(t,x)\\\q\qq -h^{(i)}(t,x, (v^j(t,x))_{j=1,m}, (\sigma^\top D_xv^i)(t,x), B_i\bar u^i(t,x))=0,\,\, (t,x)\in \esp ;\\
v^i(T,x) = g^i(x).
\end{array}  \right.
\ee
But, the family $(\bar u^i)_{i=1,m}$ belongs to $\pgc$ and solves system (\ref{secondorder-pde-inter2}). Therefore, by the uniqueness result of Proposition \ref{solvisco}, one deduces that $\bar u^i= v^i$, $\forall i=1,...,m$. 

Next we are going to show that on $[t,T]\times \Omega \times E$, $ds\otimes d\P\otimes d\l$-a.e we have: $\forall i=1,...,m$, 
\be \lb{expressionubar}\begin{array}{ll}
\bar U_{s}^{i;t,x}(e)&=
v^i(s,X^{t,x}_{s-}+\beta (s,X^{t,x}_{s-},e))-v^i(s,X^{t,x}_{s-})\\{}&=\bar u^i(s,X^{t,x}_{s-}+\beta (s,X^{t,x}_{s-},e))-\bar u^i(s,X^{t,x}_{s-}).\ea\ee
The second equality is trivial once the first one is proved. 

Note that we cannot use the result of Proposition \ref{represdesui} as we do not know whether or not the function 
$x\mapsto \bar h^{(i)}(t,x,y,z)=h^{(i)}(t,x,y,z,B_i\bar u^i(t,x))$ belongs uniformly to $\u$. However the function $(t,x)\mapsto B_i\bar u^i(t,x)$ is continuous and belongs to $\pg$, since $\bar u^i$ belongs to $\u$ and thanks to the properties (\ref{cdbeta}) and (\ref{condgamma}) on $\beta$ and $\g_i$ respectively. 

We are going to make use of the hint of Remark \ref{remimpo}. Let $(x_k)_{k\ge 1}$ be a sequence of $\r^k$ which converges to $x\in \r^k$ and let $\xkk$ and $\xk$ be the processes defined by (\ref{eq:diffusion2}) when the initial conditions are $x_k$ and $x$ respectively. Next let us consider the two following BSDEs (adaptation is w.r.t $\cF^k$):
\be \left\{ \label{mainbsde4}
    \begin{array}{l}  
\vec{\bar Y}^{k,t,x}:=(\bar Y^{i,k;t,x})_{i=1,m}\in \cS^2(\r^m),\,
{\bar Z}^{k,t,x}:=(
\bar Z^{i,k;t,x})_{i=1,m} \in \cH^2(\r^{m\times d}),\, \bar U^{k,t,x}:=(\bar U^{i,k;t,x})_{i=1,m}\in \cH^2(L_m^2(\l_k));\\\\ \forall i\in \mo,\,\,
\bar Y^{i,k;t,x}_T= g^{i}(\xk_T) \mbox{ and }\forall s\le T,\\\\
\qq d\bar Y_s^{i,k;t,x} = -h^{(i)}(s,\xk_s, \vec{\bar Y}^{k,t,x}_s,\bar Z_s^{i,k;t,x}, \int_E
\g_i(s,\xk_s,e)\{\bar u^i(s,\xk_{s-}+\beta (s,\xk_{s-},e))-\bar u^i(s,\xk_{s-})\}\l(de))ds\\\\ \qq\qq\qq +\bar Z_s^{i,k;t,x} dB_s +\int_{E} 
    \bar U_{s}^{i,k;t,x}(e)\tilde{\mu}_k(ds,de).\end{array}\right.  \ee
First by continuity and as in the proof of Step 2 of Proposition \ref{represdesui} for any $i=1,...,m$, one can check that 
$(\bar Y^{i,k;t,x},\bar Z^{i,k;t,x},\bar U^{i,k;t,x}1_{\{|e|\geq \frac{1}{k}\}})_k$ converges to 
$(\bar Y^{i;t,x},\bar Z^{i;t,x},\bar U^{i;t,x})$ in $\mathcal{S}^{2}(\R)\times \mathcal{H}^2(\mathbb{R}^{\kappa\times d}) \times \mathcal{H}^2(L^2(\l))$. Next let $((v_i^k)_{i=1,m})_{k\ge 1}$ be the sequence of continuous determinstic functions such that for any $t\leq T$ and $s\in [t,T]$,
$$\bar Y^{i,k;t,x}_s= v_i^k(s, ^k\! X^{t,x}_s)
\mbox{ and }\bar Y^{i,k;t,x_k}_s= v_i^k(s, ^k\! X^{t,x_k}_s)
, \forall i=1,...,m.
$$ 
Note that the function $v_i^k$ belongs uniformly to $\pg$, i.e. there exists a constant $C$ which does not depend on $k$ such that $|v_i^k(t,x)|\le C(1+|x|^\rho)$, 
$\forall \tx$, for some $\rho\ge 0$. On the other hand, for any $i=1,...,m$,  we have:

(i) the sequence $(v_i^k(t,x))_{k\ge 1}$ converges to $v^i(t,x)$ ; 

(ii) $\bar U^{i,k;t,x}=v^k_i(s,^k\! X^{t,x}_{s-}+\beta (s,^k\! X^{t,x}_{s-},e))-v_i^k(s,^k\! X^{t,x}_{s-})$, $ds\otimes d\P\otimes d\l_k$-ae on $[t,T]\times \Omega \times E$. 
\medskip 

\noindent Now using It\^o's formula and  the properties satisfied by $h^{(i)}$  
we obtain for some constant $C\geq 0$: 
$$\ba{ll}
&\E[|\vec{\bar Y}^{k,t,x_k}_s-\vec{\bar Y}^{k,t,x}_s|^2
+\int_s^T|\bar Z^{k,t,x_k}_r-\bar Z^{k,t,x}_r|^2ds +\int_s^T\ie|\bar U^{k,t,x_k}_r(e)-\bar U^{k,t,x}_r(e)|^2\l_k(de)]\\\\
&\qq \leq C\E[|g(^k\! X^{t,x_k}_T)-g(^k\! X^{t,x}_T)|^2]+C\E[\int_s^T
|\vec{\bar Y}^{k,t,x_k}_r-\vec{\bar Y}^{k,t,x}_r|^2 dr]\\\\ &\qq \qq +C\E[\int_s^T|^k\!X^{t,x_k}_r-^k\!X^{t,x}_r|(1+|X^{t,x_k}_r|^p+|X^{t,x_k}_r|^p)]\\ \\&\qq \qq +C\sum_{i=1,m}\E[\int_s^T|B_i\bar u^i(r,^k\!X^{t,x_k}_r)-B_i\bar u^i(r,^k\!X^{t,x}_r)|^2dr],\q \forall s\leq T.
\ea
$$
Next using Gronwall's inequality and taking $s=t$ to obtain: $\forall i=1,...,m$, 
\be \lb{estimvik}\ba{ll}
|v_i^k(t, x_k)-v_i^k(t,x)|^2&\leq 
\E[|\vec{\bar Y}^{k,t,x_k}_t-\vec{\bar Y}^{k,t,x}_t|^2]\\\\&\le 
C\E[|g(^k\! X^{t,x_k}_T)-g(^k\! X^{t,x}_T)|^2]+C\E[\int_t^T|^k\!X^{t,x_k}_r-^k\!X^{t,x}_r|(1+|X^{t,x_k}_r|^p+|X^{t,x_k}_r|^p)]\\\\&\qq\qq +C\sum_{i=1,m}\E[\int_t^T|B_i\bar u^i(r,^k\!X^{t,x_k}_r)-B_i\bar u^i(r,^k\!X^{t,x}_r)|^2dr].\ea
\ee
Finally using the estimates (\ref{estimx}) satisfied by $\xk$ and since the function $(t,x)\mapsto B_i\bar u^i(t,x)$ is continuous and belongs to $\pg$ to deduce that the right-hand side of (\ref{estimvik}) converges to $0$ as $k \rw \infty$. Henceforth the sequence 
$(v_i^k(t, x_k)-v_i^k(t,x))_k$ converges to $0$ as $k\rw \infty$ for any $i=1,...,m$. Consequently by Remark \ref{remimpo} and (i)-(ii) above we have, for any $i=1,...,m$,   
\be \lb{expressionubar2}\begin{array}{ll}
\bar U_{s}^{i;t,x}(e)=
v^i(s,X^{t,x}_{s-}+\beta (s,X^{t,x}_{s-},e))-v^i(s,X^{t,x}_{s-}), \,\,ds\otimes d\P\otimes d\l-a.e. \mbox{ in }[t,T]\times \Omega \times E.\ea\ee
which is the desired result. 
\ms

We now come back to the issue of uniqueness. Replacing in (\ref{mainbsde3}) the quantity\\ $\bar u^i(s,X^{t,x}_{s-}+\beta (s,X^{t,x}_{s-},e))-\bar u^i(s,X^{t,x}_{s-})$ with $\bar U_{s}^{i;t,x}(e)$, we deduce that the triple 
$(\vec{\bar Y}^{t,x},{\bar Z}^{t,x},\bar U^{t,x})$ verifies: $\forall i\in \mo$,
\be \left\{ \label{mainbsde4}
    \begin{array}{l} 
    \bar Y_T^{i;t,x}= g^{i}(X^{t,x}_T) \mbox{ and }\forall s\leq T,\\\\
d\bar Y_s^{i;t,x} = -h^{(i)}(s,X^{t,x}_s, \vec{\bar Y}^{t,x}_s,\bar Z_s^{i;t,x}, \int_E
\g_i(s,X^{t,x}_s,e)\bar U^{i;t,x}_s(e)\l(de))ds+\bar Z_s^{i;t,x} dB_s +\int_{E} 
    \bar U_{s}^{i;t,x}(e)\tilde{\mu}(ds,de).\end{array}\right.  \ee
It follows that $$ \forall \; i\in \mo, \quad \quad
\bar Y^{i;t,x}=Y^{i;t,x}$$ since the solution of the BSDE (\ref{mainbsde3}) is unique. 
Thus for any $i\in \mo$, $u^i=\bar u^i=v^i$ which means that the solution of (\ref{ipde-intro}) in the sense of Definition \ref{od} is unique inside the class $\u$. \qed
\section{Extensions}
A) Let us assume that for any $i\in \mo$ the functions 
$f^{(i)}$, have the following form:
$$
\forall (t,x,y,z,\z)\in \esp\times \r^{m+d}\times 
L^2(\l), f^{(i)}(t,x,y,z,\z)=h^{(i)}(t,x,y,z,\|\z\|_{L^2(\l)})$$ where 
the functions $(h^{(i)})_{i=1,m}$ are the ones defined in Section 2. Under Assumptions (H1)-(H2) on $(h^{(i)})_{i=1,m}$ and $(g^{i})_{i=1,m}$ and by Proposition \ref{existencegene} (see also Remark \ref{remexistuniq}) for any $\tx$ there exists a unique solution 
$(\vec{Y}^{t,x}, {Z}^{t,x},U^{t,x})$ of the following BSDE with jumps: 
\be \left\{ \label{mainbsdeother}
    \begin{array}{l}  
\vec{Y}^{t,x}:=(Y^{i;t,x})_{i=1,m}\in \cS^2(\r^m),\,
{Z}^{t,x}:=({Z}^{i;t,x})_{i=1,m}\in \cH^2(\r^{m\times d}), U^{t,x}:=(U^{i;t,x})_{i=1,m}\in \cH^2(L_m^2(\l));\\ \forall i\in \mo, \,
Y_T^{i}= g^{i}(X^{t,x}_T) \mbox{ and }\\
\qq dY_s^{i;t,x} = -h^{(i)}(s,X^{t,x}_s, \vec{Y}^{t,x}_s,Z_s^{i;t,x}, \|U_s^{i;t,x}\|_{L^2(\l)})ds -Z_s^{i;t,x} dB_s -\int_{E} 
    U_{s}^{i;t,x}(e)\tilde{\mu}(ds,de),\,\,\forall s\leq T.\end{array}\right.  \ee
Next by Proposition \ref{solvisco} there exist deterministic continuous functions $(\underbar u^i(t,x))_{i=1,m}$ which belong to $\pg$ such that for any $(t,x)\in \esp$, the solution of the BSDE (\ref{mainbsde}) verifies: 
\be\label{repre2}
\forall i\in \mo,\,\,\forall s\in [t,T],\,\,Y^{i;t,x}_s=\underbar u^i(s,X^{t,x}_s).
\ee
 Moreover, one can easily show that the functions $(\underbar u^i)_{i=1,m}$ belong to $\u$ and in the same way as in Section 3 the processes $U^{t,x}:=(U^{i;t,x})_{i=1,m}$ of the BSDE with jumps (\ref{mainbsdeother}) are linked to the functions $(\underbar u^i)_{i=1,m}$ by (\ref{represdesui}). Finally by the same method as in the proof of Theorem \ref{mainresult} we obtain: 

\begin{thm} \label{mainresult2}
Assume that Assumptions (H1)-(H2) are fulfilled. Then the $m$-tuple of functions $(\underbar u^i)_{i=1,m}$ defined in (\ref{repre2}) is the unique viscosity solution in the class $\u$ of the following system of IPDEs: $\forall i=1,...,m$, 
\be \left\{ \label{secondorder-pde2} \begin{array}{l} 
-\partial_t \underbar u^i(t,x)- b(t,x)^\top D_x\underbar u^i(t,x)-\frac{1}{2}\mathrm{Tr}\big(\sigma\sigma^\top (t,x) D^2_{xx}\underbar u^i(t,x)\big)-K\underbar u^i(t,x)\\\q\qq -h^{(i)}(t,x, (\underbar u^j(t,x))_{j=1,m}, (\sigma^\top D_x\underbar u^i)(t,x), B_i\underbar u^i(t,x))=0,\,\, (t,x)\in \esp ;\\
\underbar u^i(T,x) = g^i(x),
\end{array}  \right.
\ee where for any $(t,x)$, $B_i\underbar u^i(t,x)$ is given by \be\label{defbk2}\begin{array}{l} B_i\underbar u^i(t,x)=\{\int_E|\underbar u^i(t,x
+\beta(t,x,e))-\underbar u^i(t,x)|^2\l(de)\}^{\frac{1}{2}}.
\ea\ee
\end{thm}

Note that the definition of the viscosity solution of (\ref{secondorder-pde2}) is the same as the one given in Definition \ref{od} but with the new expression of $B_i\underbar u^i(t,x)$ given by (\ref{defbk2}). 

According to our best knowledge, viscosity solutions of IPDEs of type (\ref{secondorder-pde2}) have not been considered yet. \qed
\bigskip

\noindent B) In this study we have considered only standard IPDEs but our main result in Theorem \ref{mainresult} can be obtained for an IPDE, say, with one obstacle of the following type ($m=1$): 
\be \left\{ \label{obstacleipde} \begin{array}{l} 
\min\Big \{u^1(t,x)-\ell (t,x); -\partial_t u^1(t,x)- b(t,x)^\top D_xu^1(t,x)-\frac{1}{2}\mathrm{Tr}\big(\sigma\sigma^\top (t,x) D^2_{xx}u^1(t,x)\big)\\\quad \qq-Ku^1(t,x)-h^{(1)}(t,x,u^1(t,x), (\sigma^\top D_xu^1)(t,x), B_1u^1(t,x))\Big \}=0,\,\, (t,x)\in \esp ;\\
u^1(T,x) = g^1(x)
\end{array}  \right.
\ee
as far as, additionally, appropriate assumptions are assumed on the obstacle $\ell$. Mainly one should moreover suppose that $\ell$ belongs to class $\u$ and $\ell (T,x)\ge g^1(x)$. 

The general reflected BSDE with jumps associated with IPDE with obstacle (\ref{obstacleipde}), whose solution is a quadruple $({Y}^{t,x},{Z}^{t,x},{U}^{t,x},{K}^{t,x})$, is the following one: 
\begin{equation} \label{reflec-bsde} \left\{ 
\begin{array}{l} 
{Y}^{t,x}\in \cS^2(\r),\,
{Z}^{t,x}\in \cH^2(\r^d), U^{t,x}\in \cH^2(L^2(\l))  \mbox{ and }K^{t,x} \mb {continuous non-decreasing and }K_0=0\,\,;\\\\
 dY^{t,x}_s = -f(s, X^{t,x}_s, Y^{t,x}_s,Z^{t,x}_s, U^{t,x}_s)ds  -dK^{t,x}_s+ Z^{t,x}_s dB_s +\int_E U^{t,x}_s(e) \tilde{\mu}(ds,de) , \,s\leq T;\\\\
  Y^{t,x}_s \ge \ell(s,X^{t,x}_s),\,\,s\le T \;\textrm{and} \; \int_{0}^{T}(Y^{t,x}_s -\ell (s,X^{t,x}_s) )dK^{t,x}_s =0 ;\\\\
 \displaystyle{Y^{t,x}_T =g(X^{t,x}_T)}\\\\
\end{array} \right.
 \end{equation}
where $(t,x)\in \esp$ is fixed. We know that there exists a deterministic function $u^1$ which belongs to $\pgc$ such that:  $\fl \tx$,
\be \lb{reyobst}
\forall s\in [t,T], Y^{t,x}_s=u^1(s,X^{t,x}_s).
\ee
For more details one can see e.g.\cite{harrajouknine}. In the case when $\l$ is finite, the IPDE with obstacle (\ref{obstacleipde}) is already considered in \cite{Hammorlais-mai2015} without conditions (a)-(b) on $\g_1$ and $h^{(1)}$. The solution is given by $u^1$ of (\ref{reyobst}). In a forthcoming work we will deal with the case of a general L\'evy measure without assuming $\l(E)<\infty$. \qed 
\bigskip

\noindent {\bf Appendix}: Barles et al.'s definition for viscosity solution of IPDE (\ref{ipde-intro})
\ms 

\noindent In the paper by Barles et al. \cite{BarlesBuckPardoux}, the definition of the viscosity solution of the system (\ref{ipde-intro}) is given as follows.  
\begin{defn}\label{bbpdef} We say that a family of deterministic functions $u =(u^i)_{i=1,m}$, defined on $\esp$ and $\r^m$-valued and such that for any $i\in \mo$, $u^i$ is continuous, is viscosity sub-solution (resp. super-solution)
 of the IPDE (\ref{ipde-intro}) if, for any $i\in \mo$: \\
(i) $\forall x\in \r^k$, $u^i(T, x) \le g^i(x) $ (resp. $u^i(T, x) \ge g^i(x) $) ;\\
 (ii)  For any  $(t,x)\in (0,T)\times \r^k$ and any function of class $\cC^{1,2}(\esp)$ such that $(t,x)$ is a global maximum point of 
$ u^i -\phi$ (resp. a global minimum point of $ u^i -\phi$) and $(u^i-\phi)(t,x)=0$, one has
$$ -\partial_t \phi(t,x)
-\mathcal{L}^X \phi(t,x) -h^{(i)}(t,x, (u^j(t,x))_{j=1,m}, \sigma^\top(t,x)D_x\phi(t,x), B_i\phi(t,x)) \le 0
$$
(resp. 
$$-\partial_t \phi(t,x)
-\mathcal{{L}}^X \phi(t,x) -h^{(i)}(t,x, (u^j(t,x))_{j=1,m}, \sigma^\top(t,x)D_x\phi(t,x), B_i\phi(t,x)) \ge 0).
$$
The family $u =(u^i)_{i=1,m}$ is a viscosity solution of (\ref{ipde-intro}) if it is both a viscosity sub-solution and viscosity super-solution.  \qed
\end{defn} 
 
\ed 
Dear Professors

Thanks a lot for the careful reading of my paper. Mainly in this revised version, I have:

(i) Completeded the bibliography ;

(ii) Corrected the proof of Lemma 3.1

(iii) Given some details for the convergences in (3.28) which becomes (3.29) in this new version.

(iv) All the typos are corrected.

Best regards, S.Hamadene.

\end{document}

By the way, there is only still one small typos I have seen in the manuscript:
On page 8, line 5, in the definition of the norm $||\cdot||_{\delta,p}$: At the end of the expression at the right-hand side something like $]\}^{1\p}$ is missing. The correction is, of course, let at the discretion of the author..